# The S-rule and 1-d representation for the traversal of a planar graph in AEC industry

Luxin Luo[1[0000-0002-9118-273X]] , Sida Kong[2]

[1] Guiyang College of Humanities and Science, Guiyang550025, China

louise_luolx@sina.cn

[2]Guizhou University Institute of Engineering & Investigation, Guiyang550025, China

Kongsida@163.com

## Abstract

Based on two trivial observations of the AEC industry, this paper proposes a traversal method ("S-rule") and expression ("1-dimensional-graph") transformed from DFS.

This traversal method conforms to the original cognitive logic of the AEC industry, and the 1-d expression has clear language characteristics while completely retaining the topological relationship of the planar graph : a sequence of finite symbols (vocabularies) under definite rules. Moreover, the language can be restored to a standard 2-d form that is isomorphic to the original planar graph, thus ensuring its visualization characteristics.

Fragments of the 1-d language can be used as planar units for free combination and weighting, and as the data foundation to support advanced calculations including FEM and isomorphic matching. And after the 2-d graph is reduced to 1-d, any 3-d or higher-dimensional graphs can also be reduced to 1 or 2 dimensions.

The first half of this paper (Chapter 1) takes the 4X4 standard grid as an example to introduce the prototype of S-rule and 1-d expression, and gives the mapping rule from 1-d expression to its editable text form. In the second half of this paper, the rule and expression are gradually extended to non-embedded planar graphs (Chapter 2) and embedded planar graphs (Chapter 3), and the "grammar" is finally summarized (Chapter 4).

## Key words

## Background

As a special type of planar graph, in the AEC industry, the engineering drawing of the top view is mainly used for the intuitive expression of various professional components and floor spaces (especially the topological relationships and geometric features). And its “reading” rules also shape the practitioners’ cognition for the expression of engineering data. Ideally, engineering drawings can maintain 3 "degrees of freedom" in use:

1 "Free mapping": Mappings can be established between graphs of different professions and LODs(level of details) in the same area.

2 "Free division": Any area can be delineated to create a subgraph database for further editing.

3 "Free combination": Based on zoning and professional attributes, related graphs can be combined into a new one.

The purpose of this research is to establish a new type of language with both graphic and text properties for the free use of AEC planar graphs without violating the practitioners’ cognition. As the prototype of the language's grammar, the so-called "S-rule" comes from two trivial observations on AEC planar graphs:

1 Vertices and rooms can be derived from each other: the outer contour of each room can be represented as a circular vertex sequence, and each vertex can be regarded as an endpoint of the adjacent rooms’ common edge.

*(According to the AEC practice, in this study, a “face” in a planar graph is called a "room", and a subgraph containing rooms is called a "space".)*

2 Dual sequence: For general algorithms, when processing planar graph data, vertices are used as the basic unit by default, that is, a vertex sequence is generated during the process of reading the graph. For AEC practitioners, the default is to use the room object as the basic unit, that is, the process of viewing a graph will generate a room sequence.

Starting from the above observations, we put forward a hypothesis: in order to realize the unity and association of dual sequences, can we establish a special traversal order and a supporting "language"-like expression method? The basic setting of this order is: "When the last vertex of the contour vertex sequence of a room is traversed, the room is considered to have been traversed". Furthermore, the traversal order and expression method should have the following 3 characteristics:

1 Complete representation of the structure: Based on this order, an algorithm can be established to ensure that the vertex sequence can completely reproduce the room objects and their topological relationships, including multi-level embedded structures and trees.

2 "Index-free": The expression of traversal is not based on index or positioning coordinates, but is simply realized by a sequence of finite standard symbols; that is, the corresponding matching calculation between the nodes of the sequence and the nodes of the original graph is based on isomorphic reduction instead of index or coordinate checking.

3 Syntax extensibility with visual features: Based on standard symbol sequences and their reduced forms, a general interactive operation process can be established to serve the 3 degree-of-freedom requirements, and more syntaxes can be combined within professional logic.

Ultimately, this study as a whole consists of 4 parts:

1 Taking the grid as an example, the concept of “S-rule” and “1-d graph” is established, and then gradually extended to non-embedded planar graphs and embedded planar graphs, and the logic of deriving various space objects is explained.

2 The calculation principle of 1-d graph and its reduced form: starting from the vertex matching calculation, and gradually completing the calculation of the contour sequence, the reconstruction after the combination, the subgraph isomorphism and the generation of the reduced patterns.

3 Operations and calculations in interactive mode: When dealing with reduced patterns, the focus is on the establishment of general selection sets and the automatic derivation of associated subgraphs.

4 Related algorithms for 1-d graph: Including the original room calculation module as the pre-computation of the traversal process and the splitting of complex graphs, etc.

This article is the first part.

# 1 Prototypes of S-rule and 1-d Graph

## 1.1 Traversal of a grid and its one-dimensional representation

### 1 Basic Concepts of Grid and S-Rule

As a traversal order for planar graphs, the S-rule has two basic settings:

1 Once a vertex is traversed, all edges from this vertex to traversed vertices should be traversed first, and then the next vertex can be traversed.

2 The outer boundary vertices of each room are all traversed can be regarded as the room has been traversed.

According to the different levels of complexity of planar graphs, in order from simple to complex, this section first starts with the grid to explain the basic logic of the S rule. The grid is usually a node matrix composed of M rows by N columns, as shown in Figure 1.1 below, there are 3X3 rooms in the grid composed of 4X4 nodes.

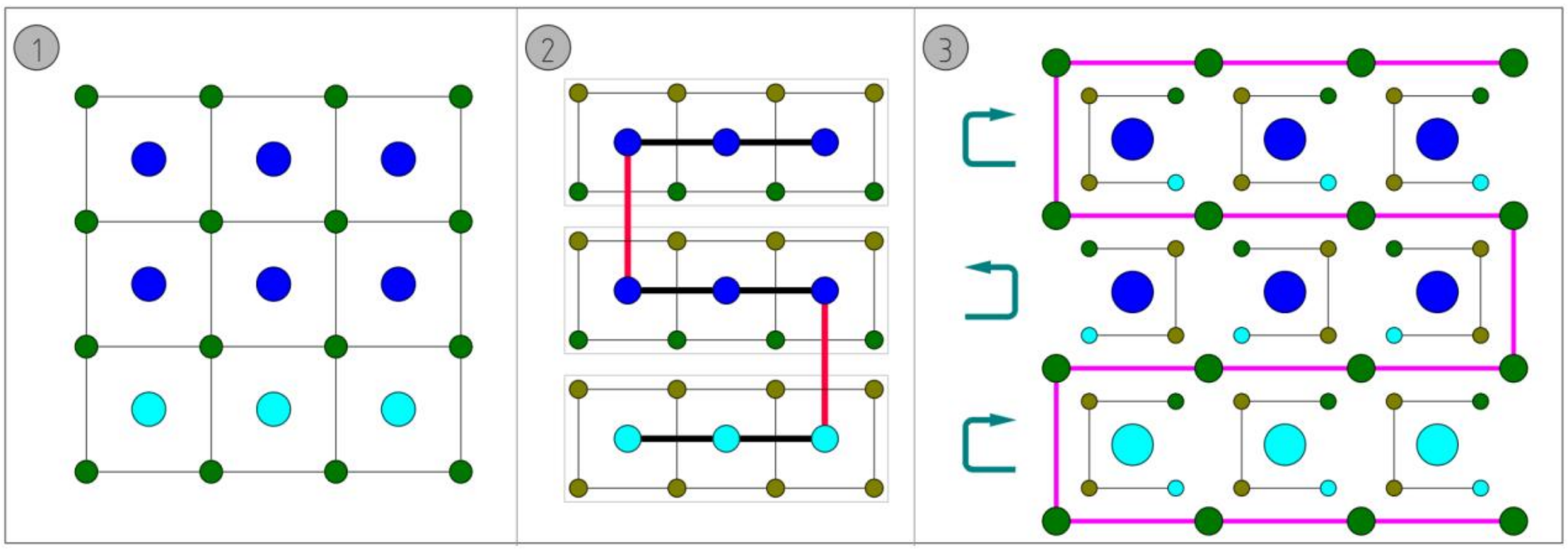

Fig.1 Strip division and traversal order of 4X4 grid

The establishment of the S-rule is based on the following 3 basic concepts:

1 Strip: Prototyped as row/column of a grid. The default "initial strip" refers to the bottom row of rooms and the two columns of vertices on the upper and lower sides of the grid, and each row of rooms together with the vertices above it is a "general strip". The subordinate rooms of each strip are combined into a "strip-subgraph", which means that the concept of strip is a division of a grid. The default traversal method of the initial strip is to start from the lower right vertex of the grid, and complete the traversal of the strip-subgraph from the outside in a clockwise manner. The S-rule states that each strip must traverse new vertices.

2 connection calculation: Refers to computing the connected subgraph of a given subgraph. Starting from the initial strip, each new strip can be obtained from the last strip based on connection calculation, which means that the division of the entire graph can be gradually obtained through it. Obviously, the strips of the grid can be obtained through initial strip and connection calculation.

3 Traversal direction: During the traversal of a planar graph, for any room, its outer boundary vertices are traversed one by one, and and finally a vertex sequence is obtained. The "traversal direction" of a room refers to the direction of the vertex sequence relative to the room itself (clockwise/anti-clockwise/null).

The concept of "traversal direction of a room" is extended to the concept of "traversal direction of a strip", that is, the S-rule adds two settings for strip and traversal direction:

1 The traversing direction of each room in each strip is consistent, and this direction is the traversing direction of the strip, while the "traversing direction" of adjacent strips is opposite.

2 2 Before a strip is traversed, it is not allowed to traverse the exclusive vertices of any new strip, but the vertices of its last strip can be traced back to derive rooms, which means that all the edges between vertices are limited to the same strip or two adjacent strips.

Based on these settings, it is possible to continue traversing vertices of the next strip from the terminal vertex of the initial level based on connection rules. Finally, a complete strip sequential traversal order can be obtained. As shown in Figure 1.2, the so-called "S" of the S-rule means that the room sequences of adjacent strips are connected end to end to form a strip chain, and the chain as a whole takes the shape of an "S". Furthermore, as shown in Figure 1.3, the vertex sequence obtained with the guidance of the room traversal direction also presents an "S" shape.

In general, the "guidance" of the strip sequence on the traversal order is dual: with the requirements of "continuity" and "reversed traversal direction" of adjacent strips, for AEC personnel, it helps to establish the room traversal order; for the algorithm, the vertex traversal order is derived based on the room traversal order.

## 2 Basic concepts of 1-d graph

In the 1-d graph system, there are 9 types of symbols used to express the underlying vertices and edges, and each symbol is expressed by a character and a number (except for the “editing-vertices”). The first letter is used to represent the basic classification, the second digit (0-9) is used to distinguish specific functions.

| Type | S/N | Class | Explanation | Symbol |
|---|---|---|---|---|
| Vertex | 1 | Position-vertex | Including general vertex and "boundary vertices" for space objects | A |
| | 2 | Room-vertex | Including general room and various virtual rooms | B |
| | 3 | Space-vertex | Including the strip and others | C |
| | 4 | call-back vertex | Used for vertex back-connection and matching to label rooms and other space objects | D |
| | 5 | vertex as edge-weight | The attribute of the edge that the vertex is attached to, does not exist in isolation | E |
| | 6 | editing symbol | Prefix/suffix marks for start, end and line breaks for 1-d graph | S |
| Edge | 1 | position edge | Including general edge, boundary edge, and various virtual edges between position vertices | F |
| | 2 | space edge | All kinds of logical edges involving space | G |
| | | | vertices in a graph | |
| | 3 | visible edge | 2-d pre-expression of 1-d graph, which is basically not involved in this article | H |

For the grid, first define the following 19 symbols for vertices and edges:

| S/N | Symbol | Explanation |
|---|---|---|
| 1 | A0 | general position vertex |
| 2 | B0 | room without any embedded subgraphs |
| 3 | C0 | strip |
| 4 | D0 | clockwise back-vertex |
| 5 | D1 | clockwise call-vertex |
| 6 | D2 | anti-clockwise back-vertex |
| 7 | D3 | anti-clockwise call-vertex |
| 8 | E0 | primary "connectivity level" of the edge to which it belongs |
| 9 | E1 | secondary "connectivity level" of the edge to which it belongs |
| 10 | F0 | direct connection between vertices at the same strip |
| 11 | F1 | direct connection between vertices at adjacent strips |
| 12 | G4 | formal edge for deriving a Cx after the end vertex of its outer boundary, or for connecting adjacent strips |
| 13 | G5 | the start mark of a Cx's constitute expression |
| 14 | G6 | the end mark of a Cx's constitute expression |
| 15 | H1 | formal edge for marking the coverage of a Bx/Cx |
| 16 | SS | start of the 1-d graph, with a "formal edge" as a suffix |
| 17 | ES | end of the 1-d graph, with a "formal edge" as a prefix |
| 18 | DS | discontinue the 1-d graph, with a "formal edge" as a prefix |
| 19 | CS | continue the 1-d graph, with a "formal edge" as a suffix |

## 3 Syntax for "1-d graph" representation of a grid

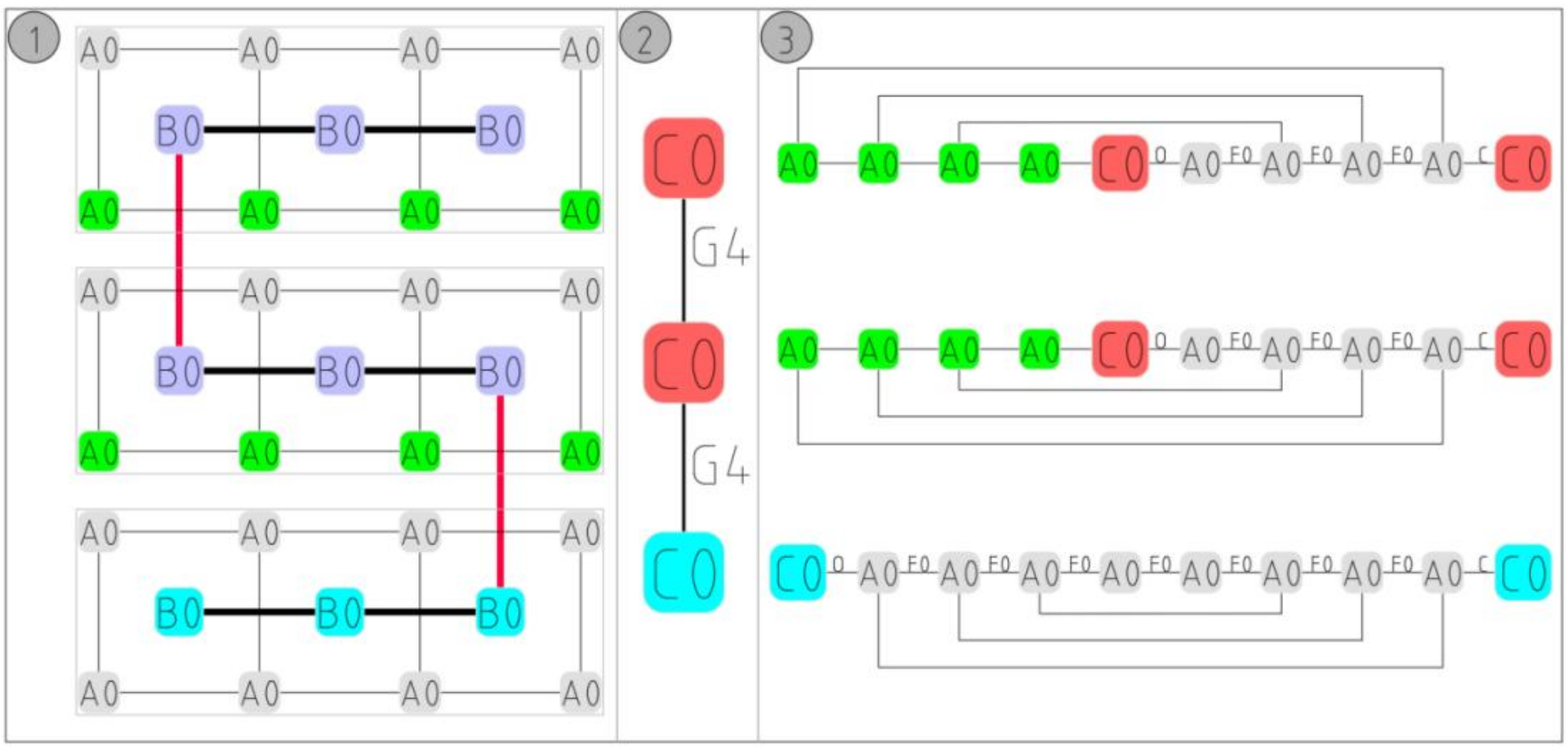

Fig.2 Strip expansion of 4X4 grid

As shown in Figure 2.1, each of the three strips has 3 rooms, the first strip has 8 new vertices, and the subsequent 2 strips have 4 traversed vertices and 4 new vertices. As shown in Figure 2.2, the three strips are arranged as a chain, connected in series with G4 symbols. As shown in Figure 2.3, the traversal of each strip (C0) is completed within a pair of C0 and a pair of "open-close" edges (O/C), and the strips are unified into a sequence composed of vertices (A0) and edges (F0). Each room corresponds to a C-shaped backtracking edge. Among them, the C-shaped edges of the last two strips are all between adjacent strips. Finally, let's arrange all 3 strips of vertex sequences into a row (1-d).

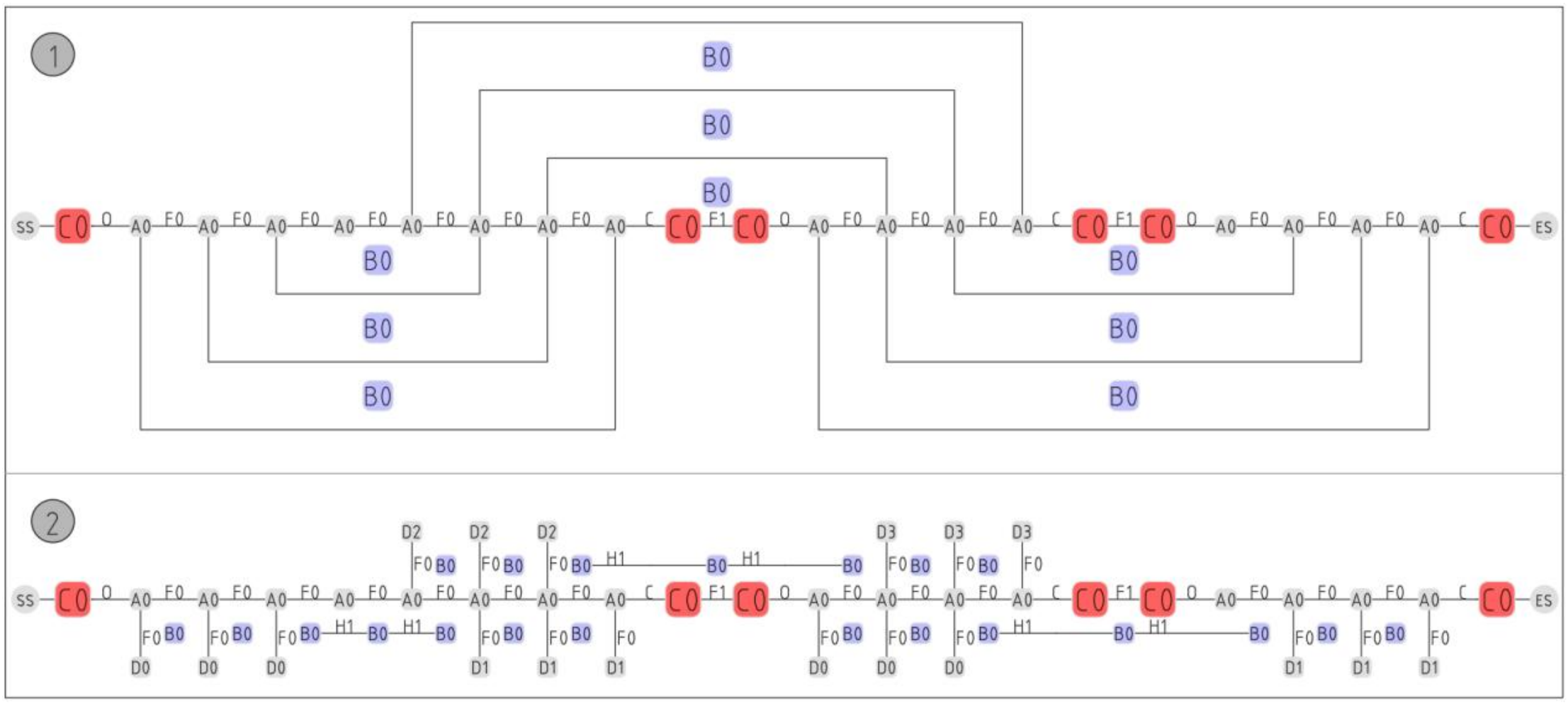

Fig.3 1-d representation of 4X4 grid

As shown in Figure 3.1, after the original data is arranged in one line, "G4" is replaced by "F1" to refer to the relationship between vertices at the boundary of adjacent strips. This form is called "restored format of 1-d graph". Furthermore, as shown in Figure 3.2, each backtracking edge is expressed by a pair of leading edges and vertices. The rooms (B0) are also placed in the corresponding position of the edges, and are attached with the scope edges of the rooms (H1). The resulting form is the "standard format for 1-d graph", and the "grammar" of this format can be briefly summarized as follows:

1 The horizontal central sequence is composed of Ax and Cx vertices, and connected by Fx and Gx;

2 The macro format of each Cx expansion is fixed: "Cx + Open" at the beginning and "Close+Cx" at the end;

3 Each room (C-shaped edge) is expressed by a call-back pair and its accompanying "vertical edge pair". The clockwise room is represented using a pair of D0/D1 vertices and a pair of "downward" identical edges. The anti-clockwise room is represented using a pair of D2/D3 vertices and a pair of "upward" identical edges.

Obviously, the restored format of 1-d graph is isomorphic to the original grid, whereas the standard format and the restoration format are equivalent in fact (that is, they can be transformed into each other based on simple calculations). The main significance of the standard format is to reduce the C-shaped edges to 1-dimension, which in turn enables it to express more complex planar structures.

## 1.2 Simplification, Extension and Textualization of 1-d Graph

### 1 Four general simplifications

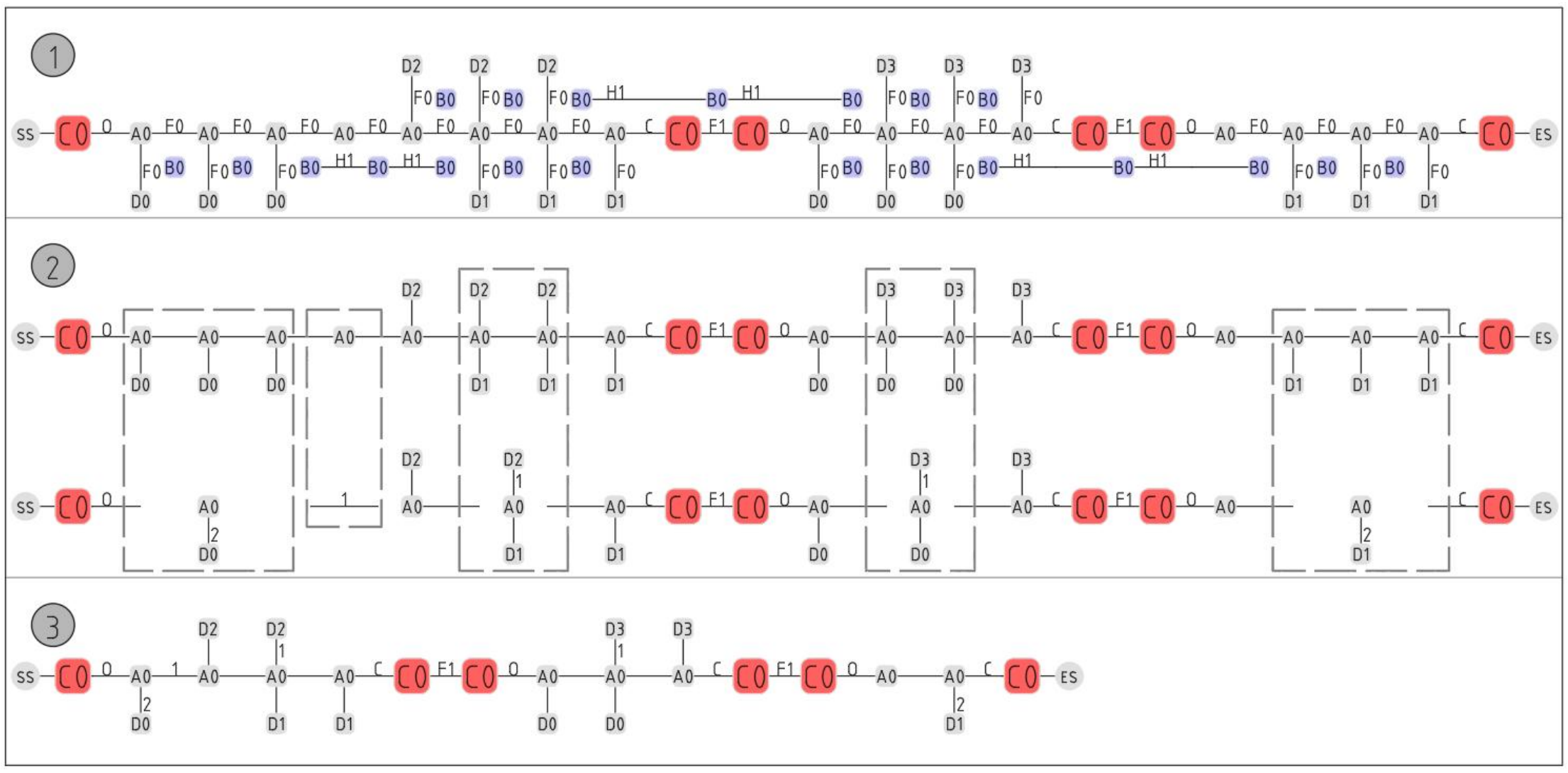

Fig.4 Simplification of 1-d graph

As shown in Figure 4, four types of information in the standard format can be simplified in stages:

1 The "F0" characters between the position-vertices and the "H1" characters between the room-vertices can be omitted, and only the lines are kept;

2 The room label (Bx) and the association between rooms (Hx) can be omitted as a whole;

3 For a continuous "...F0-A0-F0-A0-F0..." sequence without accompanying call-back vertices, it can be simplified as an edge with a value, the specific value is the number of A0 vertices contained in the sequence;

4 There may be continuous A0 sequences with identical call-backs, and the horizontal/vertical edges are all F0. The sequence can be simplified as an A0 vertex with its call-back, and the value of simplified symbols is assigned to its vertical edge. When both upward and downward vertical edges exist, the value is assigned to the upward one.

The final result is shown in Figure 4.3, which is a "simplified format".

## 2 Extension to the grammar

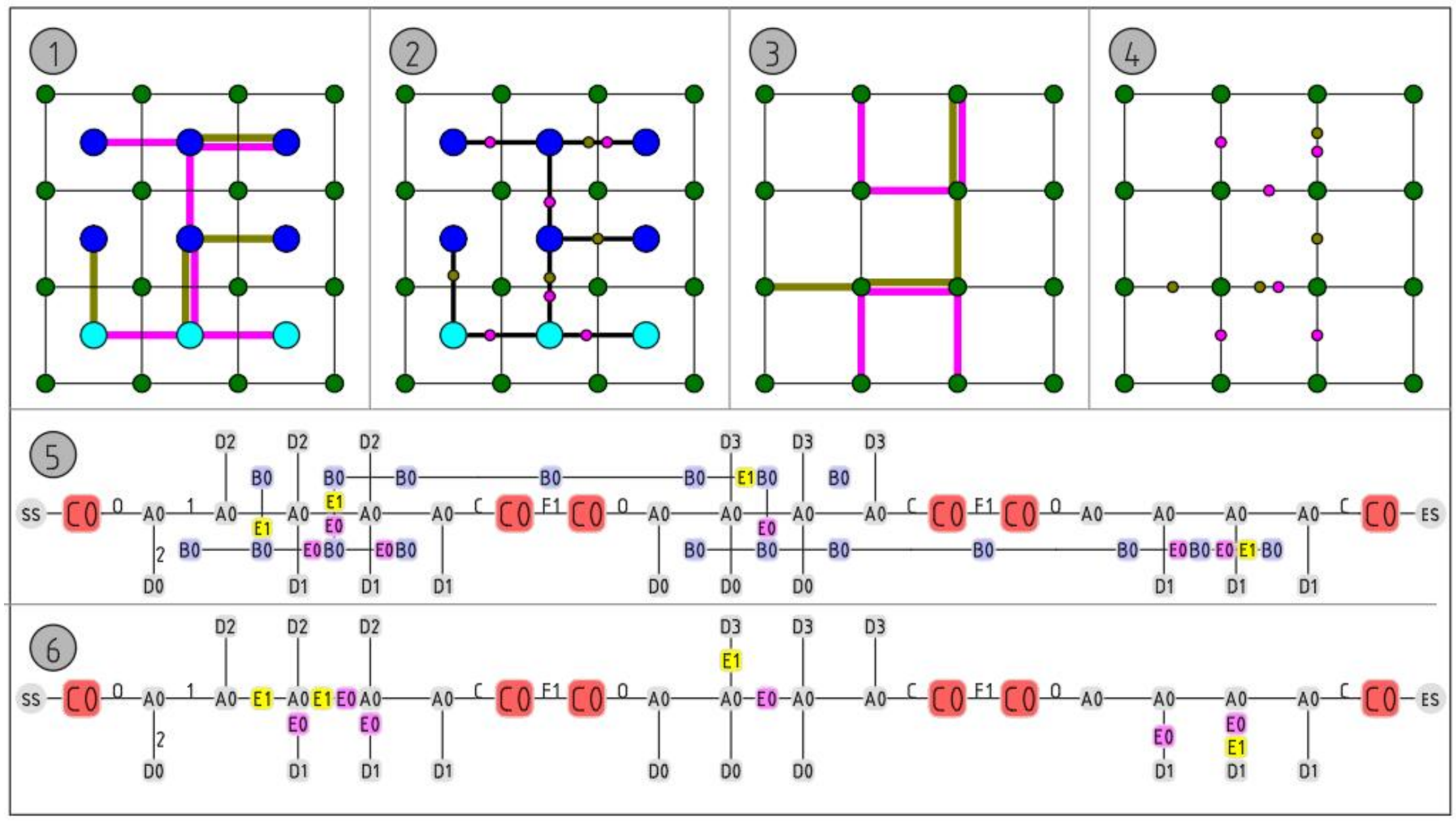

Fig.5 1-d representation of room connectivity or edge weighting

There may be specific communication information between adjacent rooms, as shown in Figure 5.1. Assuming that there are two kinds of connectivity information, such information can be equivalently converted into the weighting of the edges between adjacent rooms, as shown in Figure 5.2.

For the "simplified format", the attached object of information can be transformed into the public edges of the rooms, as shown in Figures 5.3 and 5.4.

The final results are shown in Figures 5.4 and 5.6, where the symbols (E0 and E1) are directly attached to the corresponding edges in the 1-d graph. For vertical edge pairs, the "syntax" is set to attach the Ex to the "call" side, and this setting will also affect the simplified results.

## 3 Textualization rules for 1-d graph

"Byte unit" means to regard the "main vertex" (Ax/Cx) as the base point, together with its left edge and all symbols above and below it as a unit package, including call-backs, vertical edges, derived symbols (Bx/Cx), derived associations (Hx) and assignments (Ex) on each edge. The "textualization" of a 1-d graph refers to converting it as a sequence of byte units into equivalent and reversible text strings one by one, and performing textual editing including "wordwrap".

The information composition of the "string" corresponding to a "byte unit": prefix, the backbone, 2 "simplification" values, derivative part, suffix, no more than 6 parts:

1 Prefix: SS or CS, for the start of an entire sequence or a single line;

2 Backbone: No more than 3 items, followed by main edge and Ax/Cx, downward edge and Dx, upward edge and Dx;

3-4 "Simplification" values: the values attached to the horizontal edge and the vertical edge;

5 Derivative part: no more than 3 items, followed by the lower left edge and Bx/Cx, the upper left edge and Bx/Cx, and the middle edge;

6 Suffix: ES or DS, for the end of an entire sequence or a single line;

Obviously, except for the two "simplification" values, the length of the other segments is 2*N. And 1-d graph system sets 4 basic rules for "string" mapping:

1 The pair of edge and vertex: The strings are arranged in order of edge, Ex-sequence and Ax/Cx/Dx. If the association is F0, it can be omitted.

2 The delimiter between the derivative part and the simplification value is ".", and the delimiter between adjacent "byte units" is "_".

3 If there is no "derived part" in the "byte unit", the string can omit this segment. If the "derivative part" only has the upper left subsection, then the lower left subsection is replaced by "XX", while other derived missing items are not expressed.

4 If the two "simplification" values are both 0, both segments can be omitted. However, if one item is not 0, the two segments need to be fully expressed.

With a line of 7 "byte units" as the standard length, the cases in Figure 5.5 and Figure 5.6 are line-wrapped and string-mapped, and the results are shown in Figure 6 and Figure 7:

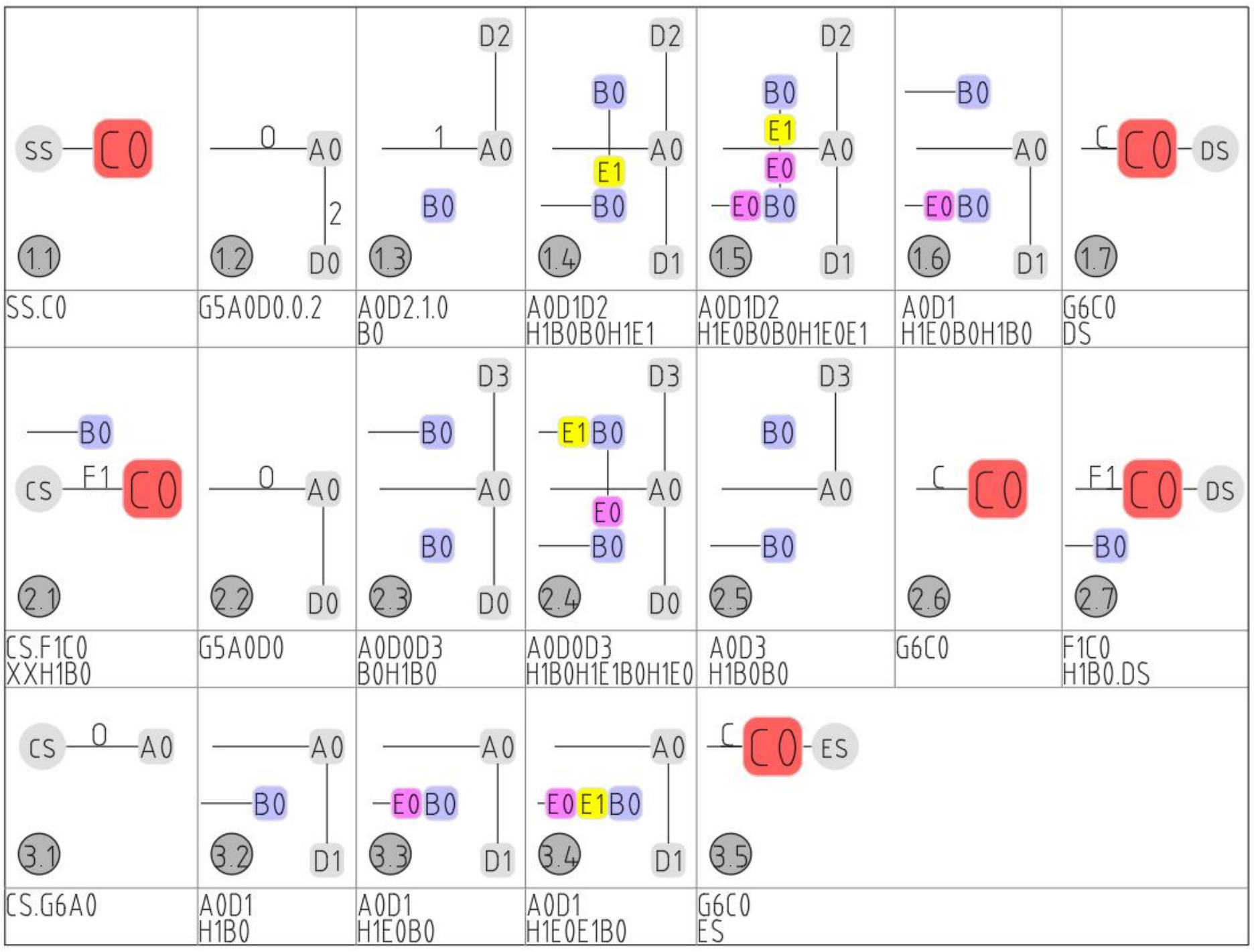

Fig.6 Byte units of 1-d graph with derived parts

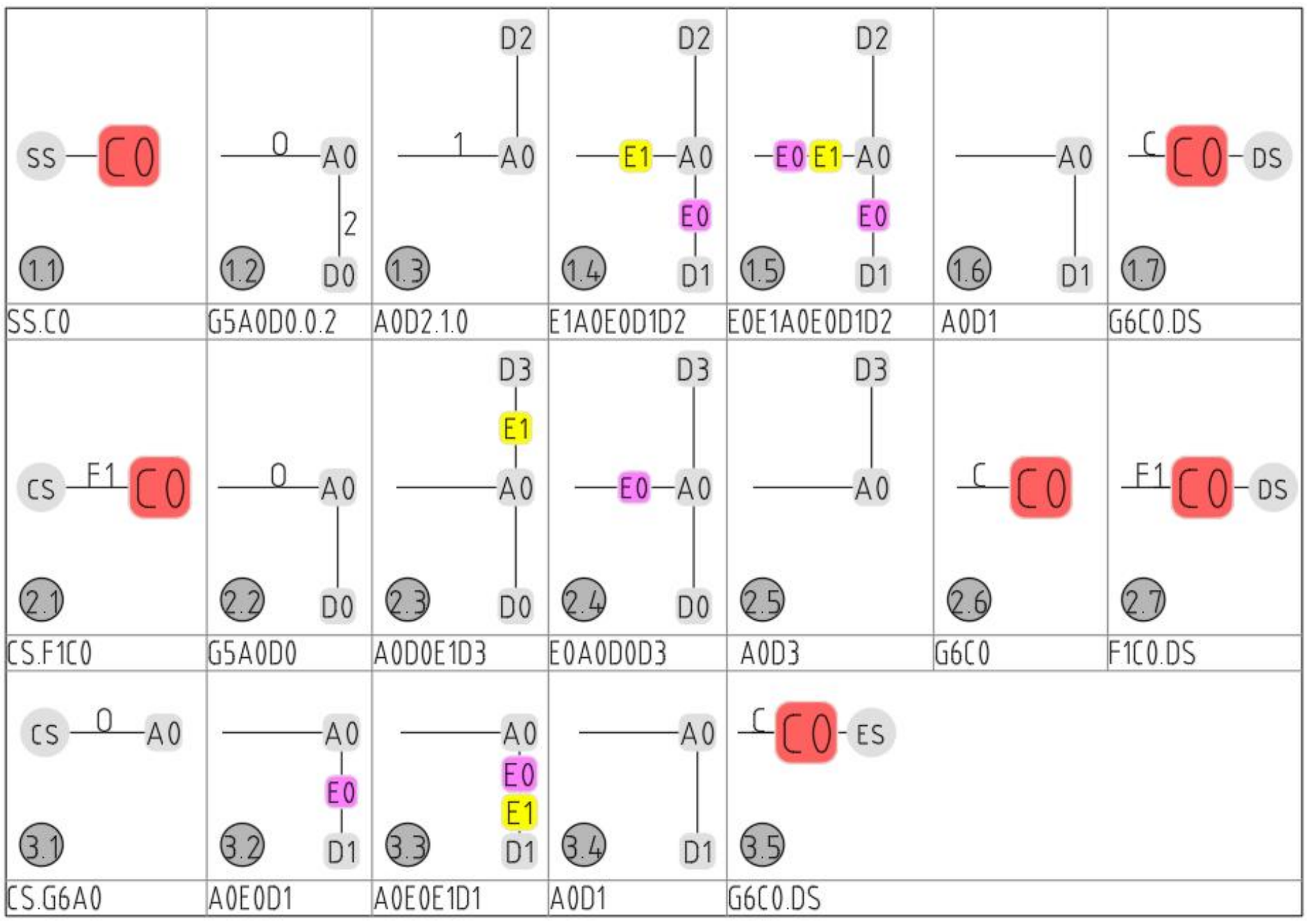

Fig.7 Byte units of 1-d graph without derived parts

The final expression of strings:

| Class | textual expression |
|---|---|
| derived format | SS.C0_G5A0D0.0.2_A0D2.1.0.B0_A0D1D2.H1B0B0H1E1_<br>A0D1D2.H1E0B0B0H1E0E1_A0D1.H1E0B0H1B0_G6C0.DS |
| | CS.F1C0.XXH1B0_G5A0D0_A0D0D3.B0H1B0_A0D0D3.H1B0H1E1B0H1E0_ |
| | A0D3.H1B0B0_G6C0_F1C0.H1B0.DS |
| | CS.G5A0_A0D1.H1B0_A0D1.H1E0B0_A0D1.H1E0E1B0_G6C0.ES |
| non-derived format | SS.C0_G5A0D0.0.2_A0D2.1.0_E1A0E0D1D2_E0E1A0E0D1D2_A0D1_G6C0.DS |
| | CS.F1C0_G5A0D0_A0D0E1D3_E0A0D0D3_A0D3_G6C0_F1C0.DS |
| | CS.G5A0_A0E0D1_A0E0E1D1_A0D1_G6C0.ES |

# 2 Grid to E0-graph

**Chapter introduction**

There are 5 sections in this chapter. Section 1 introduces the tree structure, Sections 2-4 introduce the combined-grid, and Section 5 introduces the E0-graph. In the order of chapters, space concepts including tree (C7), branch (C1), main route, inward-subgraph (C4), C4-vertex (C2), C4-chain (C3) and other corresponding symbols will be introduced in turn.

Each section consists of 3 subsections: Subsection 1 introduces the core definition, specific case and its composition under the S-rule; Subsection 2 introduces specific traversal rules and supporting concepts; Subsection 3 introduces the symbols corresponding to concepts and the 1-d representation of the case.

## 2.1 Tree

**1 Concept and Case**

Obviously, a "tree" has two typical characteristics:

A Connectivity, that is, based on association calculation, the entire graph can be traversed from any node.

B There is no room, that is, the path between any two nodes is unique.

The specific case is shown in Figure 8.1:

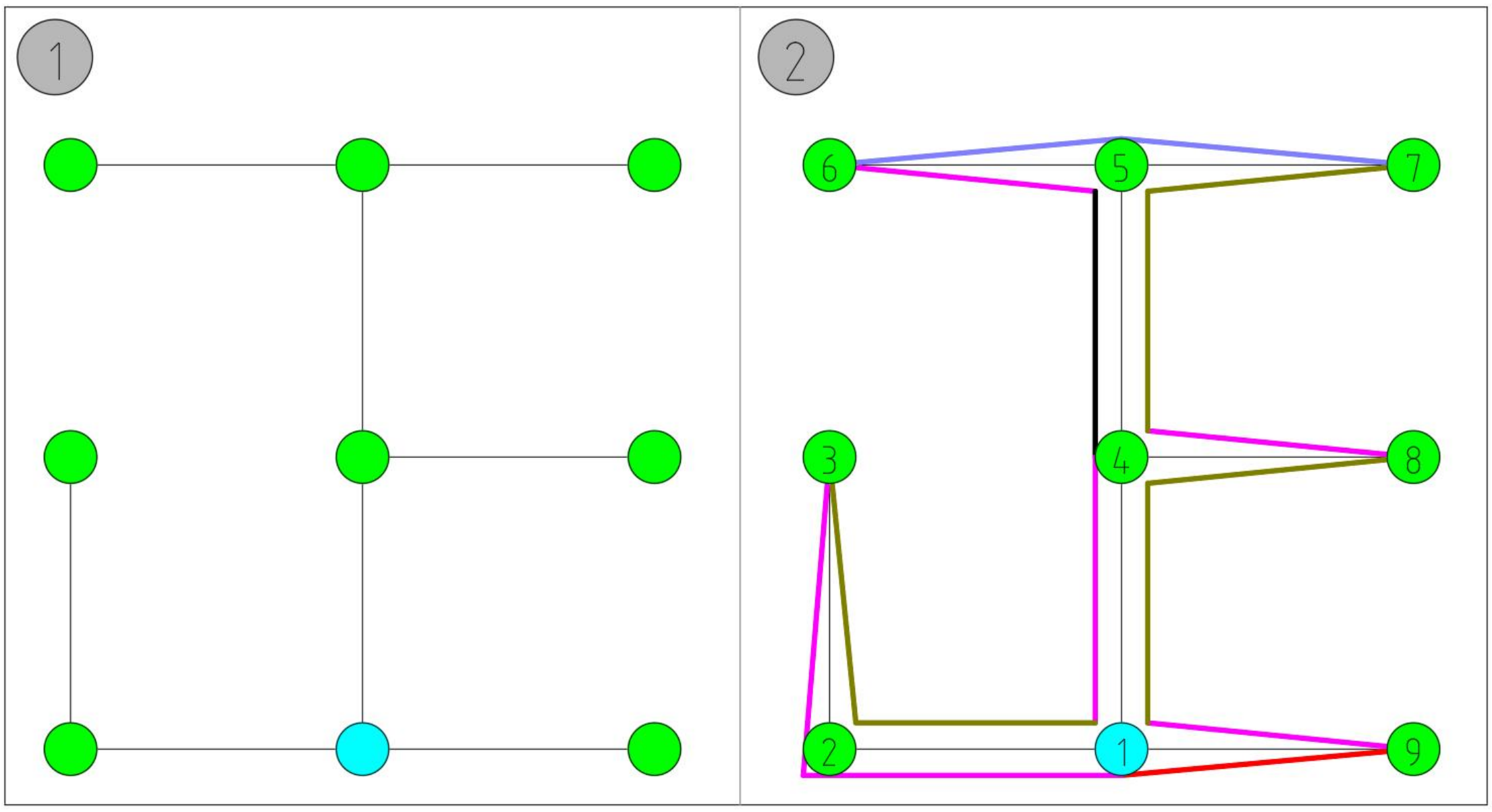

Fig.8 DFS traversal of a tree

**2 backward-edge, virtual-room, identical-edge and sibling-edge**

The S-rule follows the DFS for trees. The initial data is a vertex, and the traversal starts from the vertex, along the traversal-direction of its strip, and finally returns to the initial vertex after visiting the entire tree. During the traversal process, no new strips or other space objects are set, and the "boundary" of a tree is the vertex sequence obtained by S-rule.

Backward-edge, virtual room and identical-edge: As shown in Figure 8.2, vertex 3 returns to vertex 1. Under the DFS, after a leaf vertex is traversed, it needs to go back to the nearest traversed non-leaf vertex. And the S-rule sets that each backward edge derives a virtual room, which is marked on the 1-d graph by attaching a call-back pair. Since the vertices at both sides of the call-back are identical, the vertical edge attached to the call-back will use the new F3 symbol.

Sibling-edge: As shown in Figure 8.2, vertex 6 is associated with vertex 7. When the backward target of a leaf-vertex is directly associated with itself, sibling-edge is used to simplify the expression and directly associate to the forward vertex. And because no call-back matching is required, sibling-edges do not derive virtual rooms.

**3 Syntax extensions and case expressions for 1-d graph**

For trees, the 1-d graph system adds the following 5 symbols:

| S/N | symbol | Explanation |
| --- | --- | --- |
| 1 | C7 | Tree |
| 2 | G8 | Backward-edge |
| 3 | B1 | Virtual-room |
| 4 | F3 | Identical-edge |
| 5 | F4 | Sibling-edge |

The 1-d representation is shown in Figure 9, including restored format, standard format and simplified format. Obviously, since vertex 1 as the start needs to bear the byte occupancy of multiple back vertices, it also needs to use F3 as horizontal edges to complete its own chain expansion.

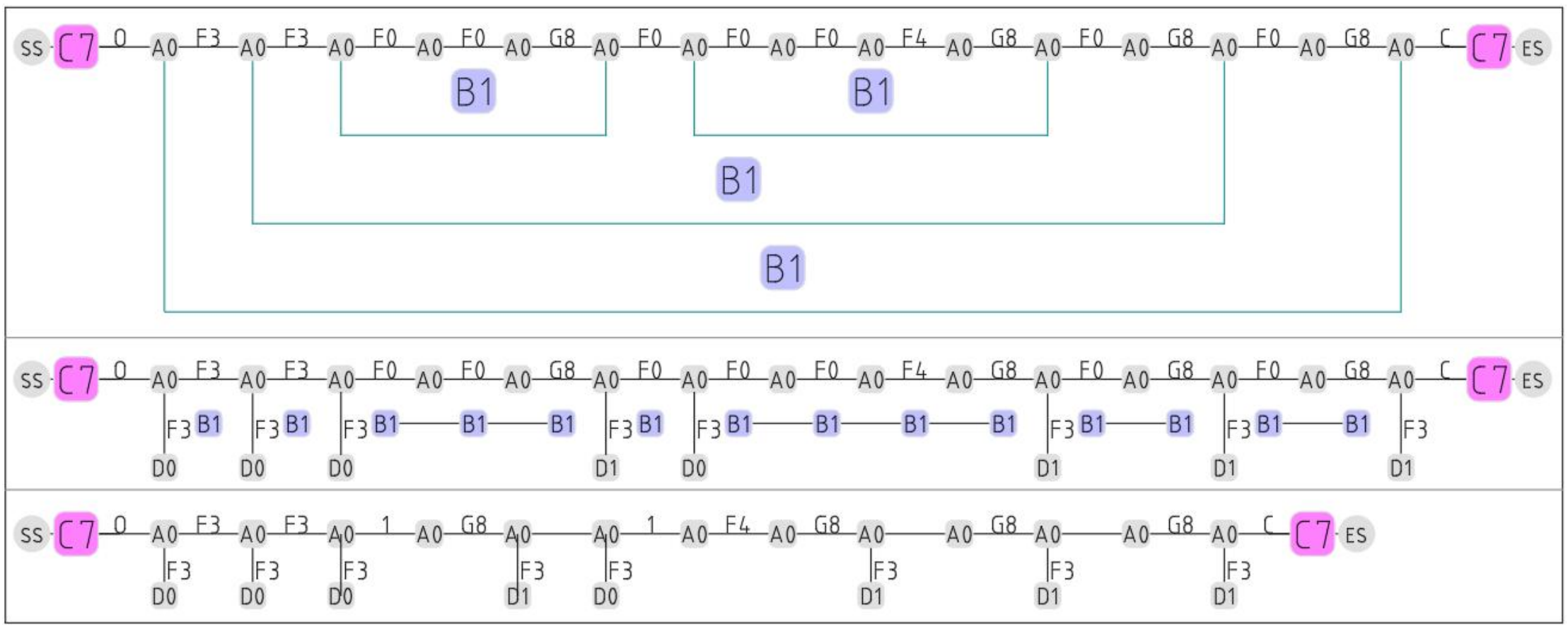

Fig.9 1-d representation of a tree

## 2.2 combined-grid

### 1 Concept and Case

The combined grid refers to a planar graph composed of grids, which has two typical characteristics:

A "Room continuity": There is a continuous path between any two rooms, and the adjacent rooms on the path have common edges;

B "Combination" will not enclose non-standard rooms, that is, all rooms have exactly 4 vertices.

The default initial data of the combined grid is a row of rooms on the outside, as shown in Figure 10.1:

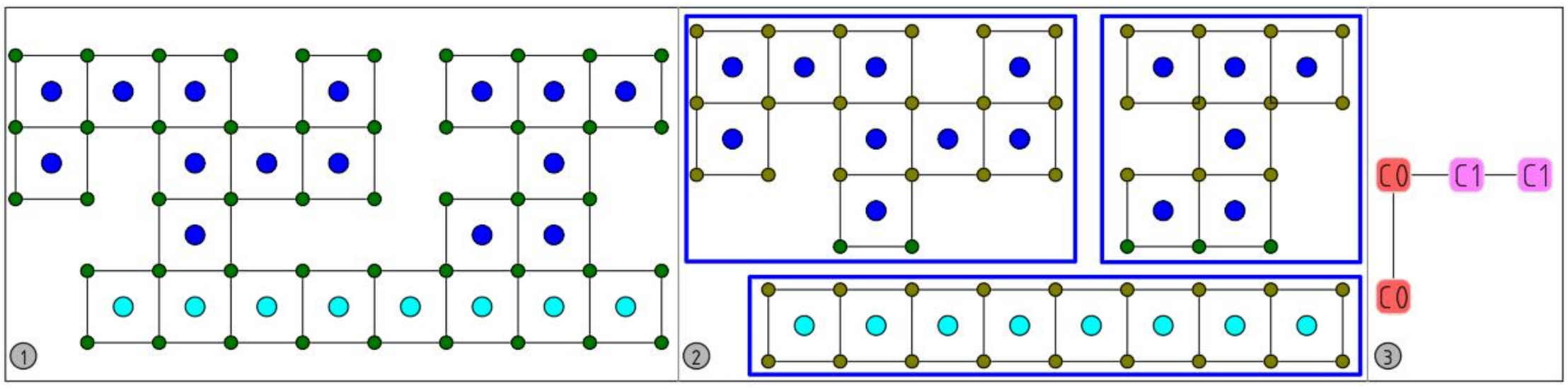

Fig.10 Preliminary division of a combined-grid

**2 branch and virtual-edge**

Branch: Refers to the sub-grid derived in the process of strip-division. As shown in Figure 10.2, a strip obtained through connection calculation is discontinuous; then, as shown in Figure 10.3, the strip derives 2 branches (symbol C1).

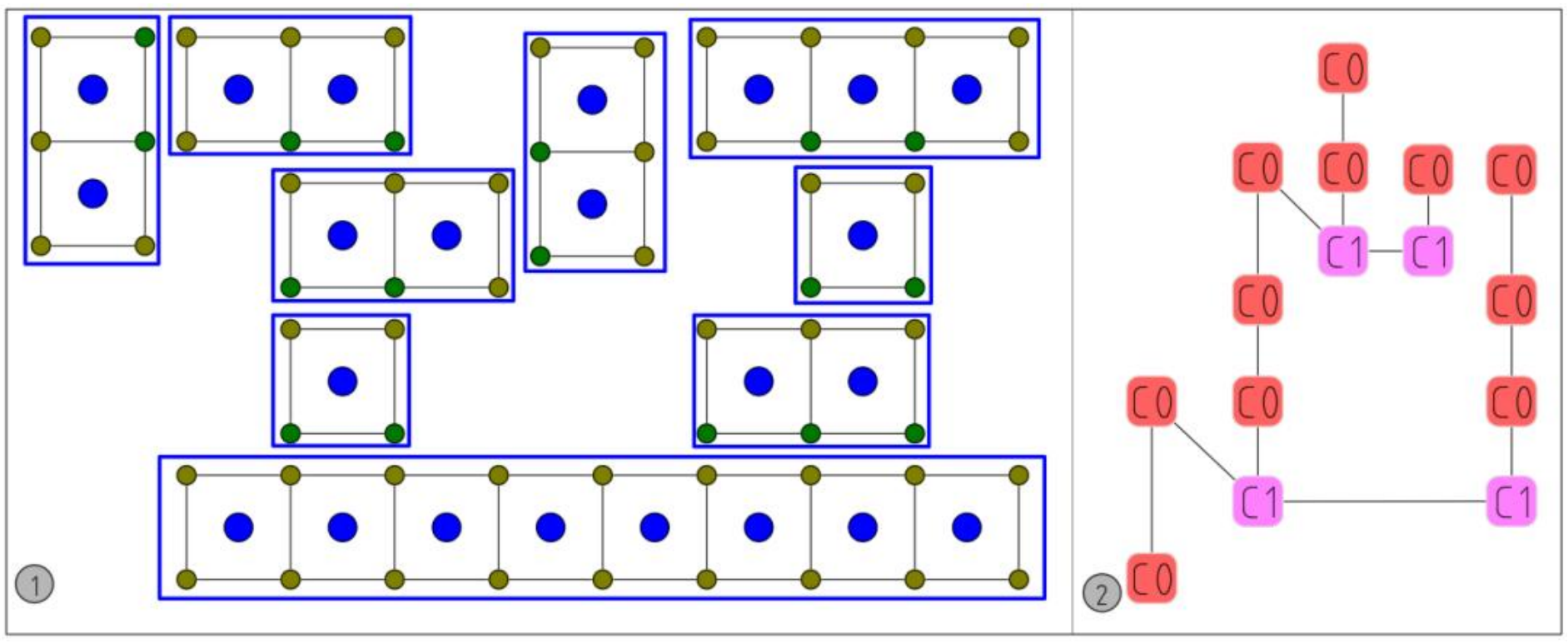

Fig.11 Complete division of a combined-grid

Furthermore, as shown in Figure 11.1, within each branch, the strip-division is continued, and the initial data in each branch is the sub-grid composed of continuous room chains connectes with the last strip. Obviously, "strip" and "branch" are derived from each other, and the division result of a combined-grid is a tree structure composed of strips and branches, as shown in Figure 11.2.

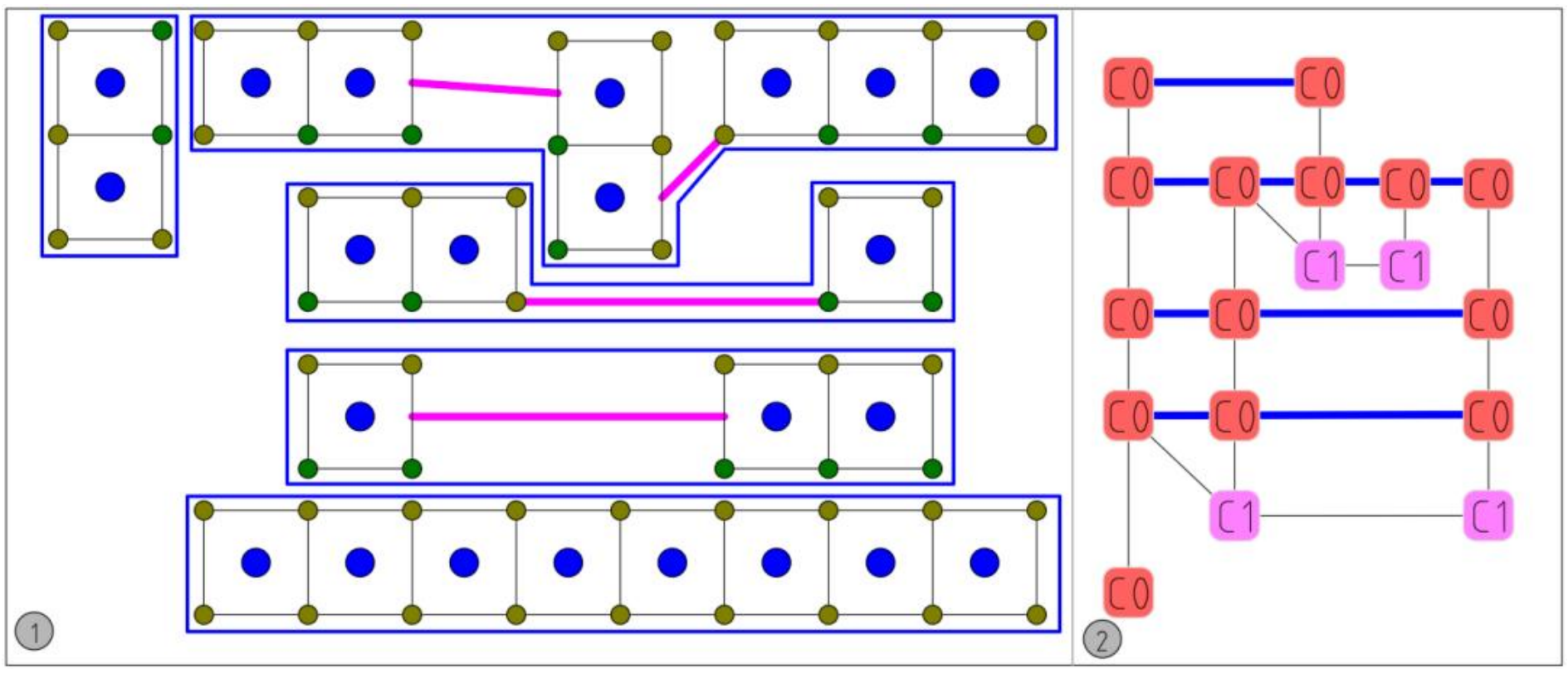

Fig.12 Virtual edges and virtual rooms of a combined-grid

Virtual-edge: When the overall setting is to avoid the appearance of branches, the virtual-edges are automatically added to the vertex sequence of the strip to generate virtual rooms, so as to ensure the formal continuity of the room sequence. Finally, as shown in Figure 12, the virtual-edges ensures that the division result is still only a strip sequence.

## 3 Syntax extensions and case expressions for 1-d graph

For combined-grids, the 1-d graph system adds the following 7 symbols:

| S/N | Symbol | Explanation |
|---|---|---|
| 1 | C1 | Branch |
| 2 | G2 | Backward edge for connecting to the next strip |
| 3 | G3 | Backward edge for connecting to the next branch(C1) |
| 4 | B3 | Virtual room for marking a G3 call-back pair |
| 5 | G0 | Virtual-edge for maintaining a strip's continuity |
| 6 | B2 | Virtual room for marking a G0 call-back pair |

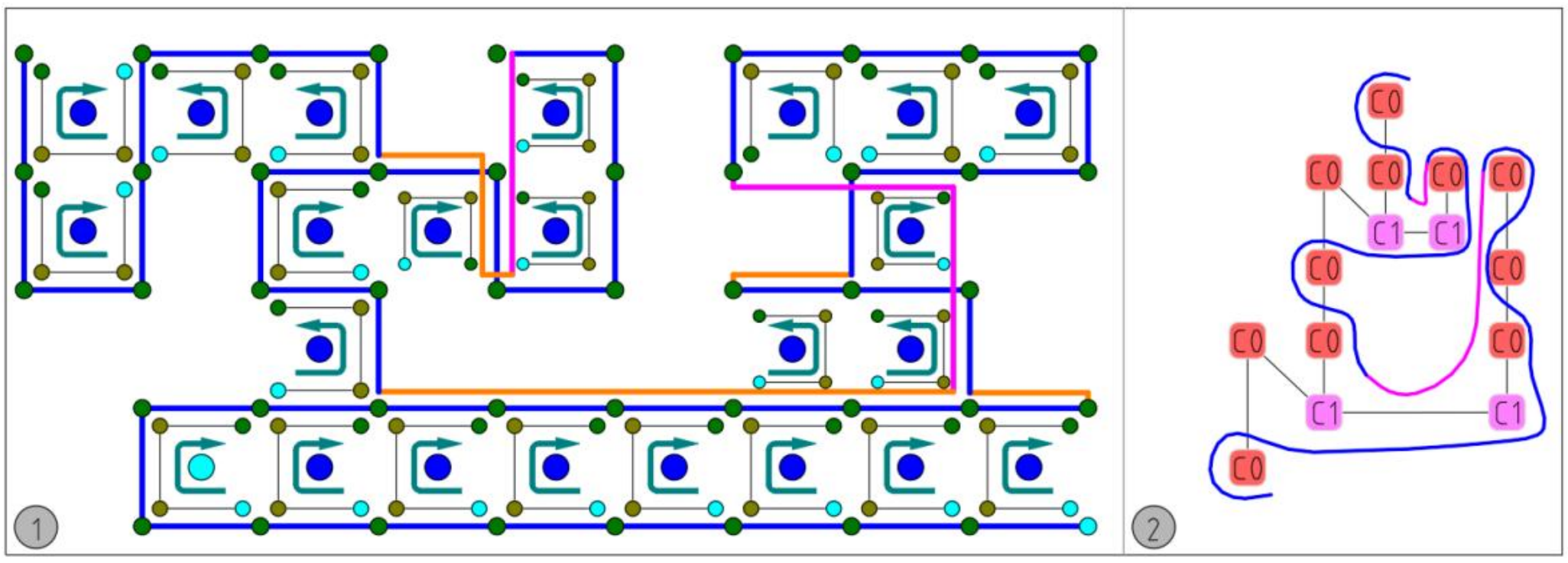

Fig.13 Traversal order of a combined-grids in branch mode

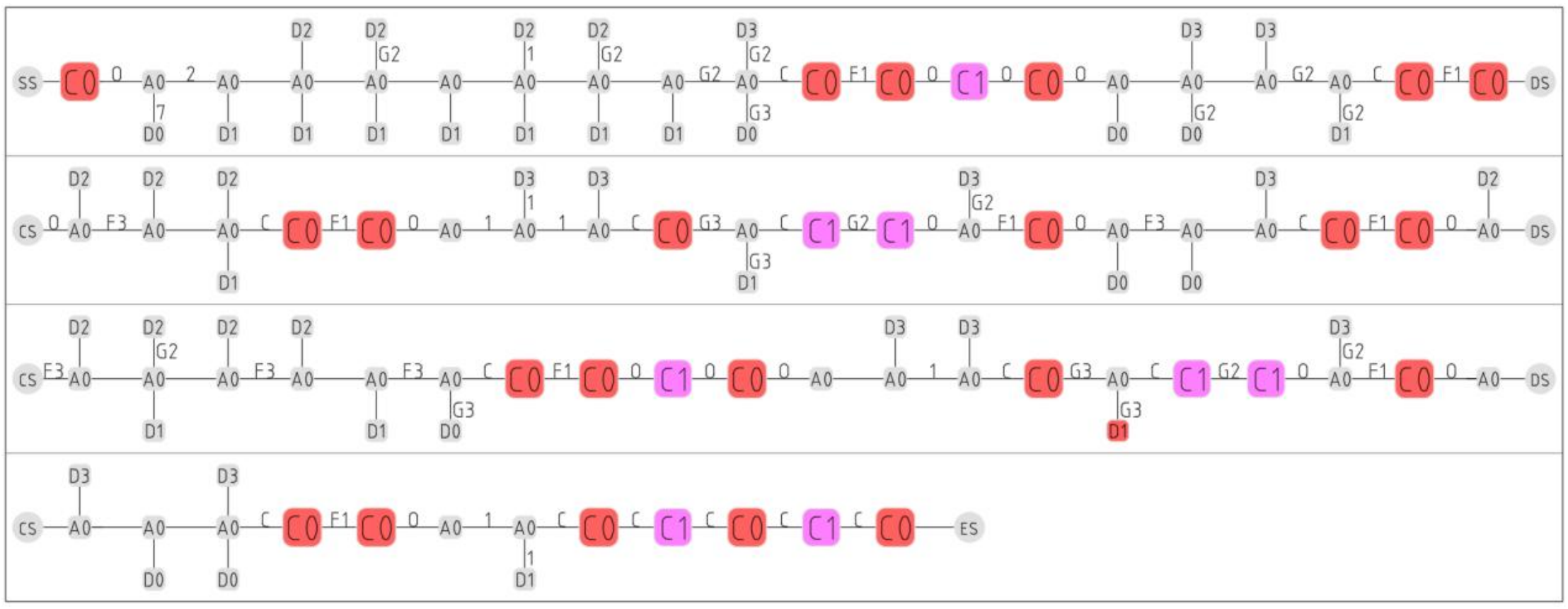

Fig.14 1-d representation of branch mode

Figure 13.1 shows the order of vertex traversal in the branch mode, while the internal logic is shown in Figure 13.2, which traverses the tree composed of strips and branches according to DFS. The 1-d representation is shown in Figure 14.

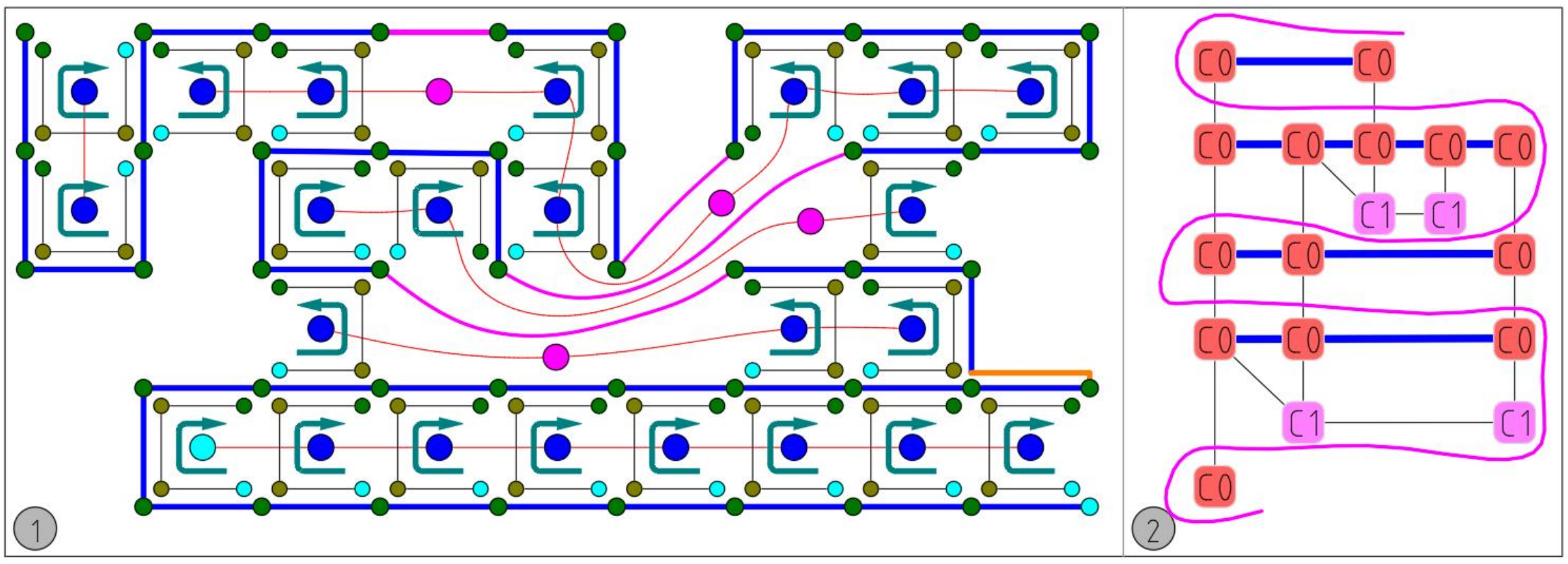

Fig.15 Traversal order of a combined-grids in virtual-edge mode

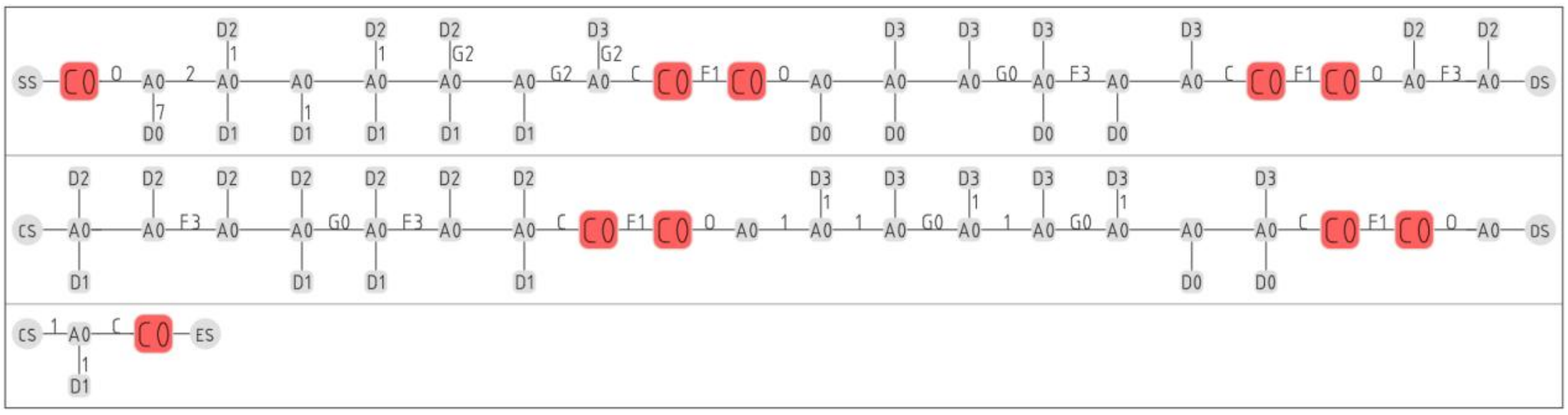

Fig.16 1-d representation of virtual-edge mode

Figure 15.1 shows the vertex traversal order in the virtual-edge mode, while the internal logic is shown in Figure 15.2, traversing the strip chain in an S-shape. The 1-d representation is shown in Figure 16.

## 2.3 outward principle

### 1 Concept and Case

For grids and combined-grids, the outward principle refers to the rule that assumes the initial data starts from an interior subgraph and then derives outward. Taking the 4X4 grid as an example, the whole is divided into two strips, as shown in Figures 17.1 and 17.2:

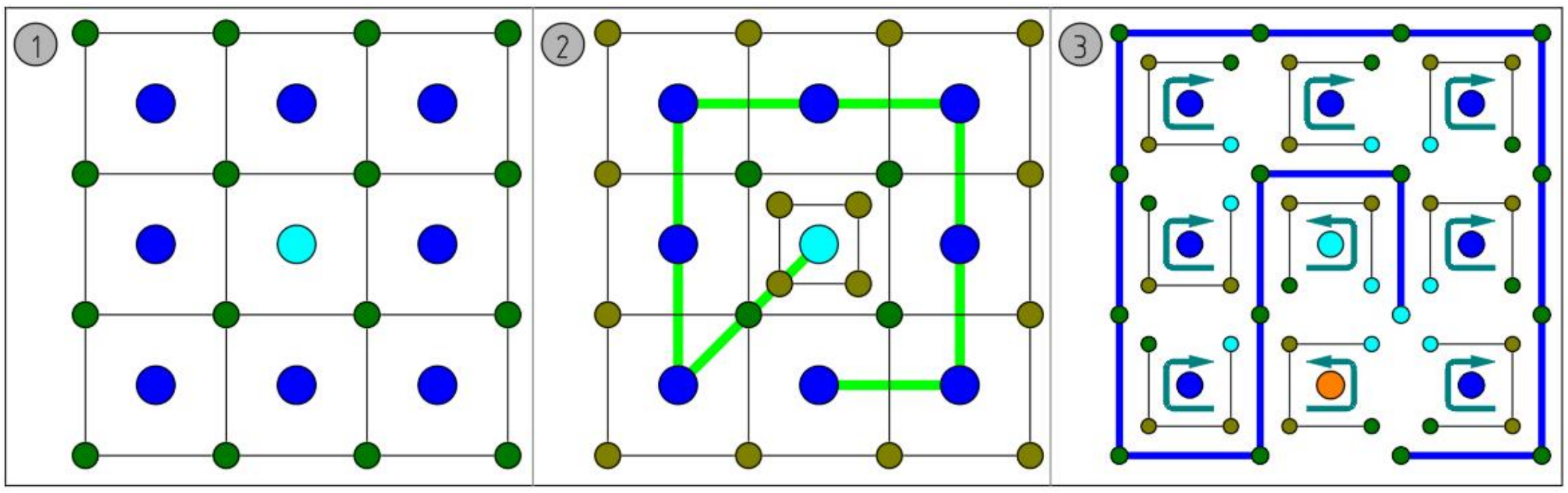

Fig.17 Strip division and traversal order of an outward grid

**2 Reversed room at the end of an outward circular strip**

As shown in Figure 17.3, for an outward circular strip, when the end room is generated, the traversal direction is opposite to the strip, which is called a reversed room. In this case, the 1-d graph is set to be marked with a reversed call-back pair, as shown in Figure 18.2, which is the corresponding restored format.

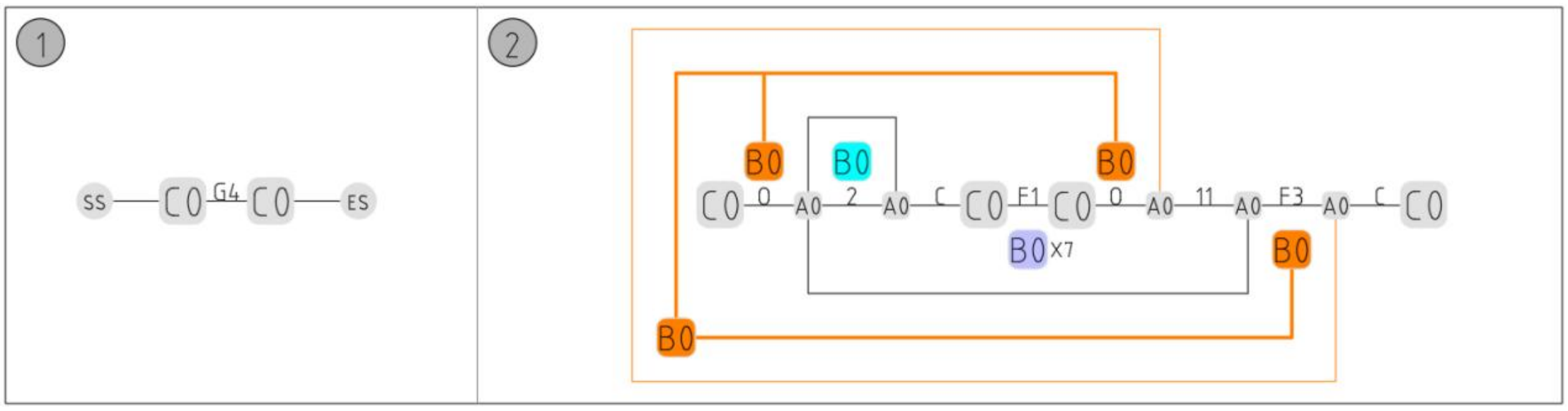

Fig.18 Logical expression of the reversed room in an outward grid

**3 Syntax extensions and case expressions for 1-d graph**

For reversed rooms, the 1-d graph system adds the following 4 symbols:

| S/N | Symbol | Explanation |
|---|---|---|
| 1 | D5 | reversed clockwise call-vertex |
| 2 | D6 | reversed clockwise back-vertex |
| 3 | D7 | reversed anti-clockwise call-vertex |
| 4 | D4 | reversed anti-clockwise back-vertex |

The 1-d representation of the outward case is shown in Figure 19:

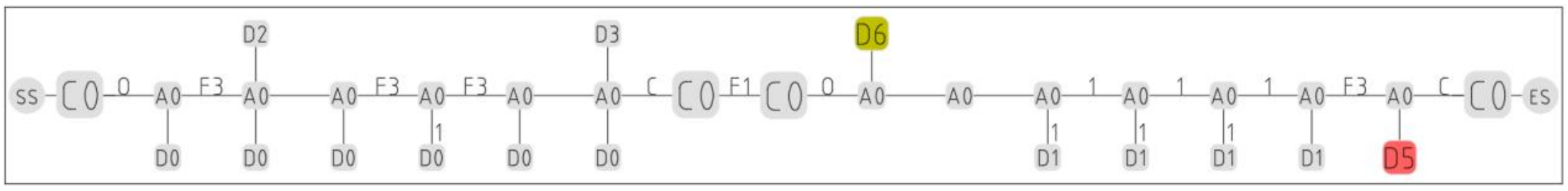

Fig.19 1-d representation of outward grid

## 2.4 inward principle

**1 Concept and Case**

For grids and combined-grids, the inward principle refers to the rule of traversing the outer boundary first, and then traversing strip by strip from the outside to the inside. The case is shown in Figure 20.1.

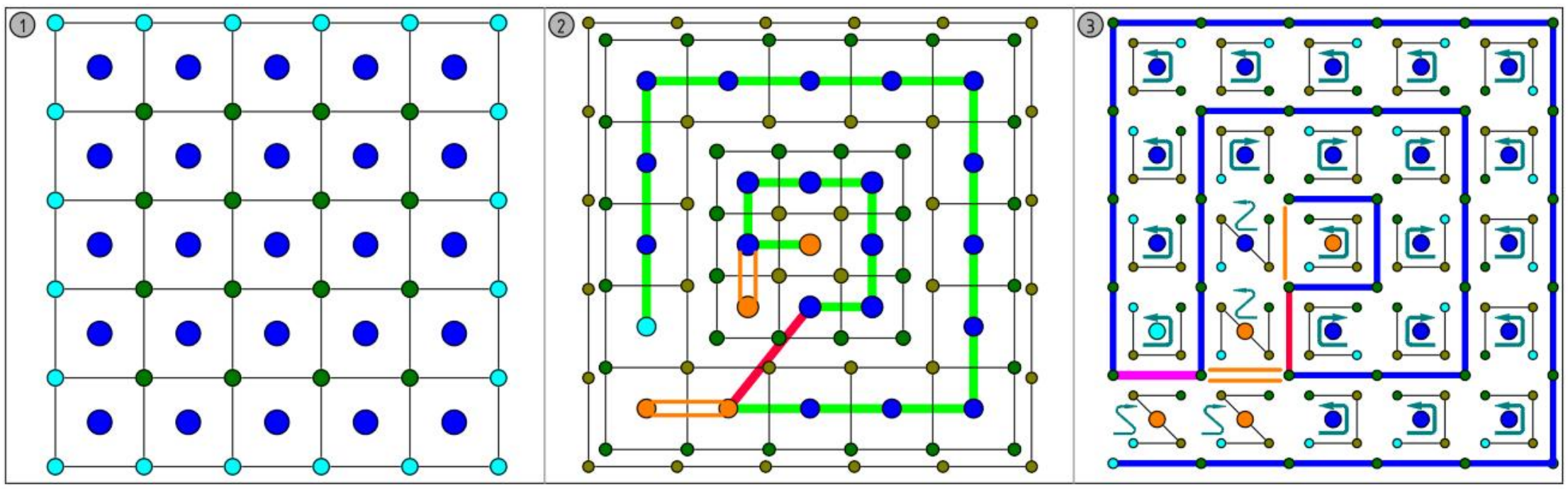

Fig.20 Strip division and traversal order of an inward 6X6 grid

**2 Boundary of inward-subgraph, S-room and inward-tree**

Boundary of inward-subgraph: For a combined-grid, the S-rule sets that when it is traversed as an inward-subgraph, its outer boundary is traversed first. That is to say, the S-rule regards it as a "big room", and the boundary-vertex-sequence is a preliminary expression. And when continuing to traverse its internal vertices, the boundary is also the first half of the first strip.

S-room: The orange room-markers shown in Figure 20.2 and Figure 20.3 include the S-rooms and the reversed room. During the inward traversal process, because the S-rule sets the "strip-traversal-direction" to remain unified, some vertices do not visit all edges in time during traversal, which means that some "call-back pairs" are placed behind. Finally, the traversal direction of the relevant rooms are S-shaped.

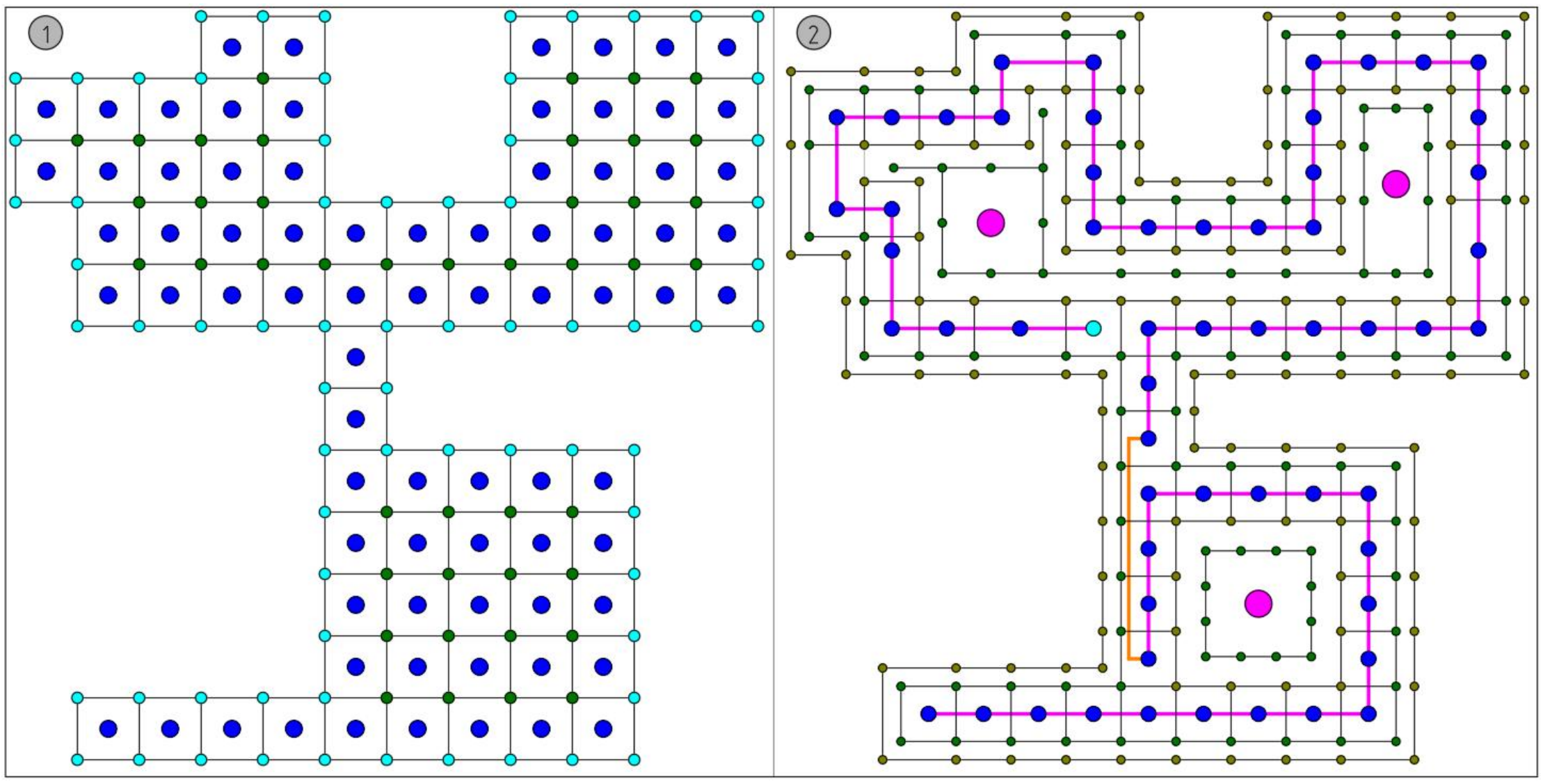

Fig.21 Strip division of an inward combined-grid

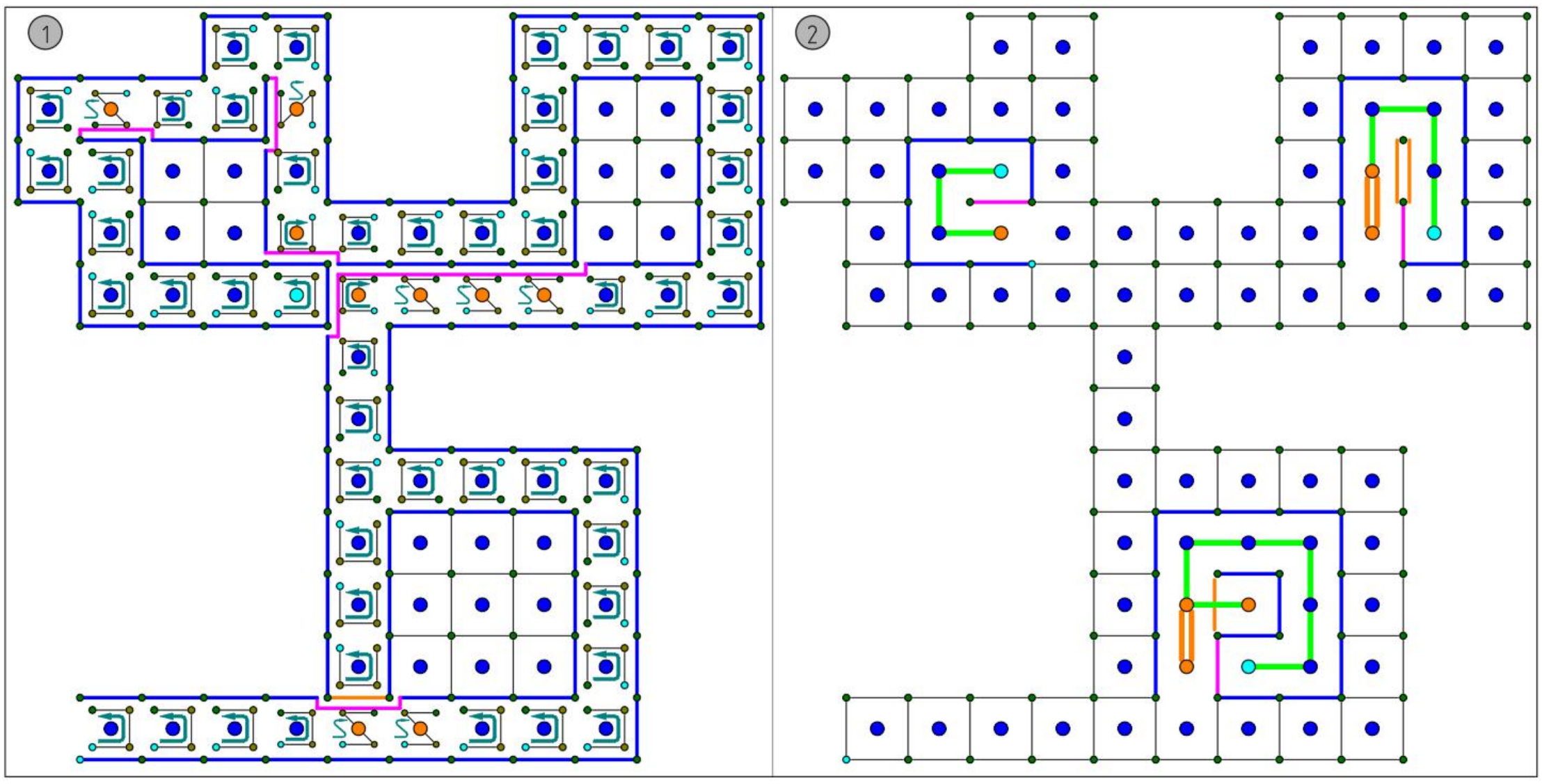

Fig.22 Traversal order of an inward combined-grid

Inward-tree: As shown in Figure 21.2, when the combined grid derives a new strip inwardly, the strip may also be discontinuous. And the boundarys of the secondary inward-subgraphs have been derived in the second half of the last strip. This means that the last strip has derived multiple secondary inward-subgraphs. Furthermore, the secondary inward-subgraphs also follow the "inward principle" when traversing. Therefore, it is obvious that, similar to the concept of branch, the strip and inward-subgraphs are derived from each other to form a tree structure. The vertex traversal order of the first strip is shown in Figure 22.1, and the vertex traversal of the three inward-subgraphs is shown in Figure 22.2.

## 3 Syntax extensions and case expressions for 1-d graph

For inward-subgraphs, the 1-d graph system adds the following 3 symbols:

| S/N | Symbol | Explanation |
|---|---|---|
| 1 | C4 | Inward-subgraph |
| 2 | B4 | Boundary of an inward-subgraph |
| 3 | F2 | the edge used to connect the inner and outer boundary of a C4's first strip |

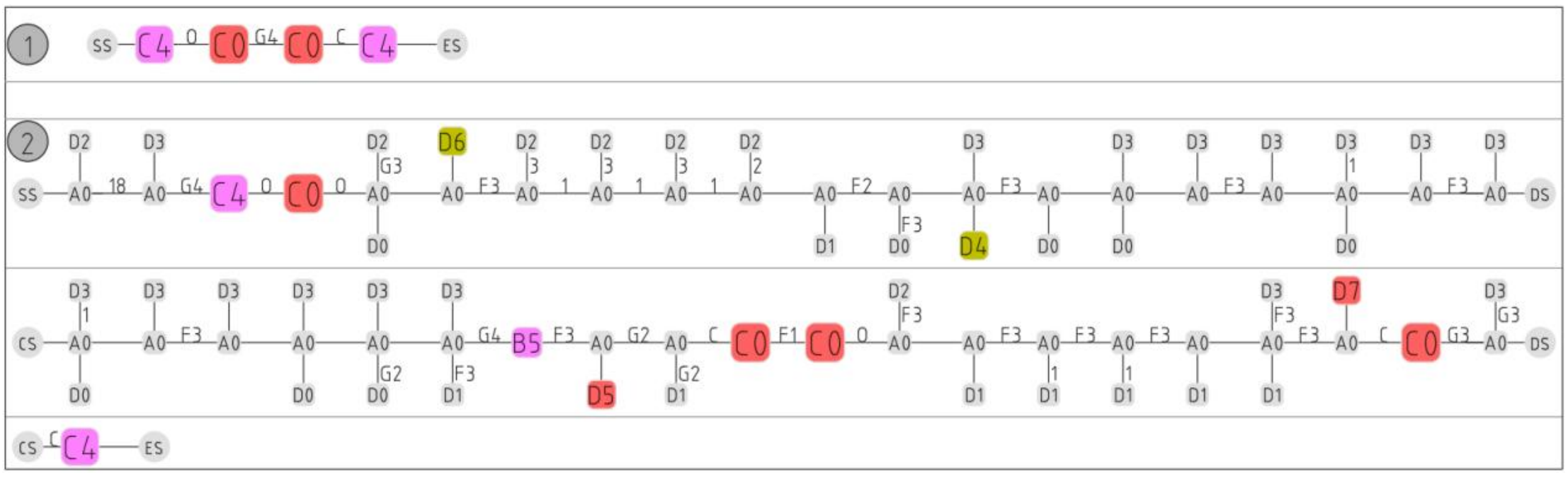

Fig.23 1-d representation of inward 6X6 grid

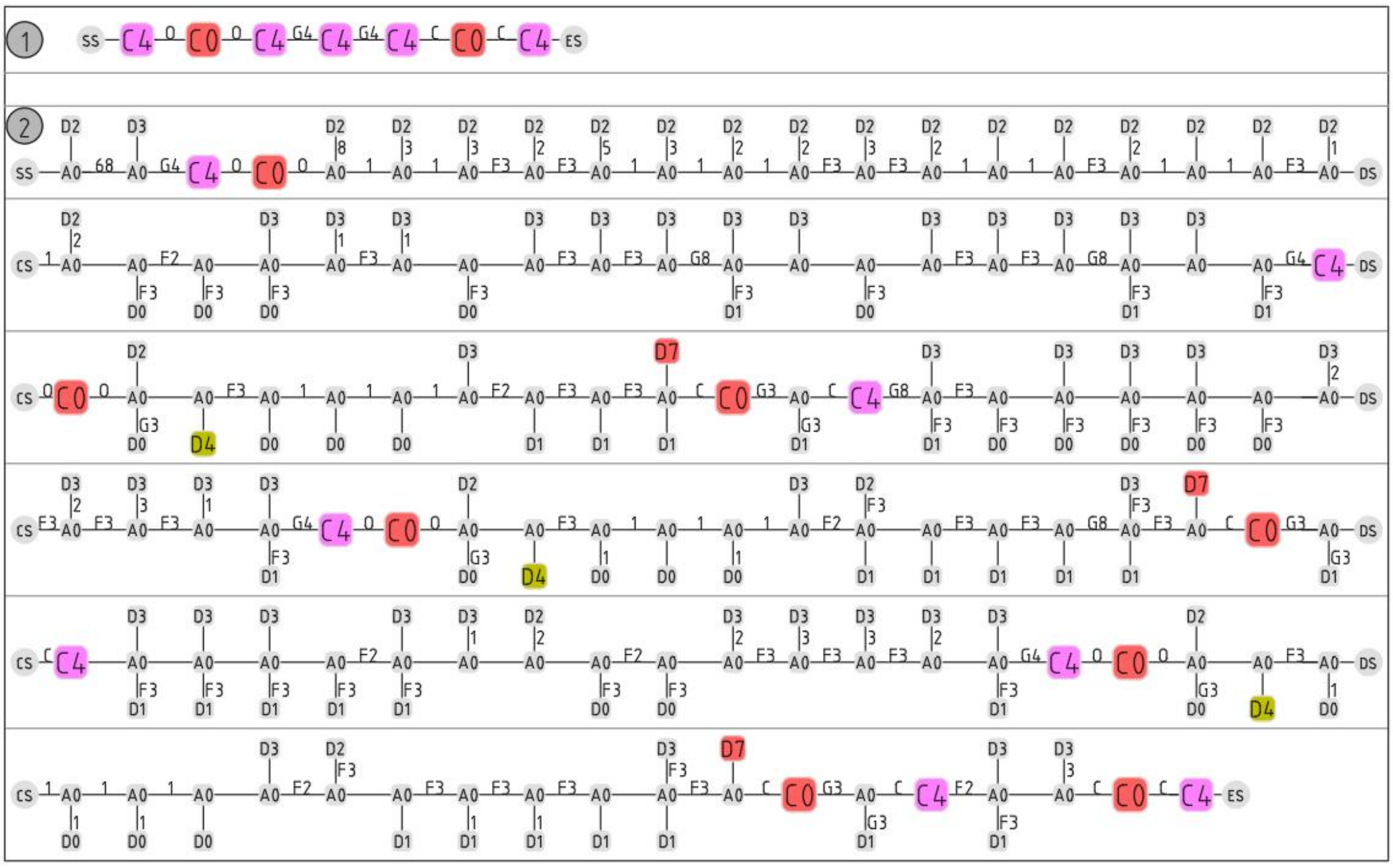

Fig.24 1-d representation of inward combined-grid

As shown in Figure 23.1 and Figure 24.1, the overall space composition is very simple, while the unfolded 1-d representation shown in Figure 23.2 and Figure 24.2 is relatively complex.

## 2.5 E0-graph

### 1 Concept and Case

As shown in Figure 25.1, the E0-graph has 3 typical characteristics:

A "Room Continuity", same as required for combined-grid.

B "No embedment", which means that any position-vertex is not inside a room..

C No external tree-edges, which means that all edges are room boundaries.

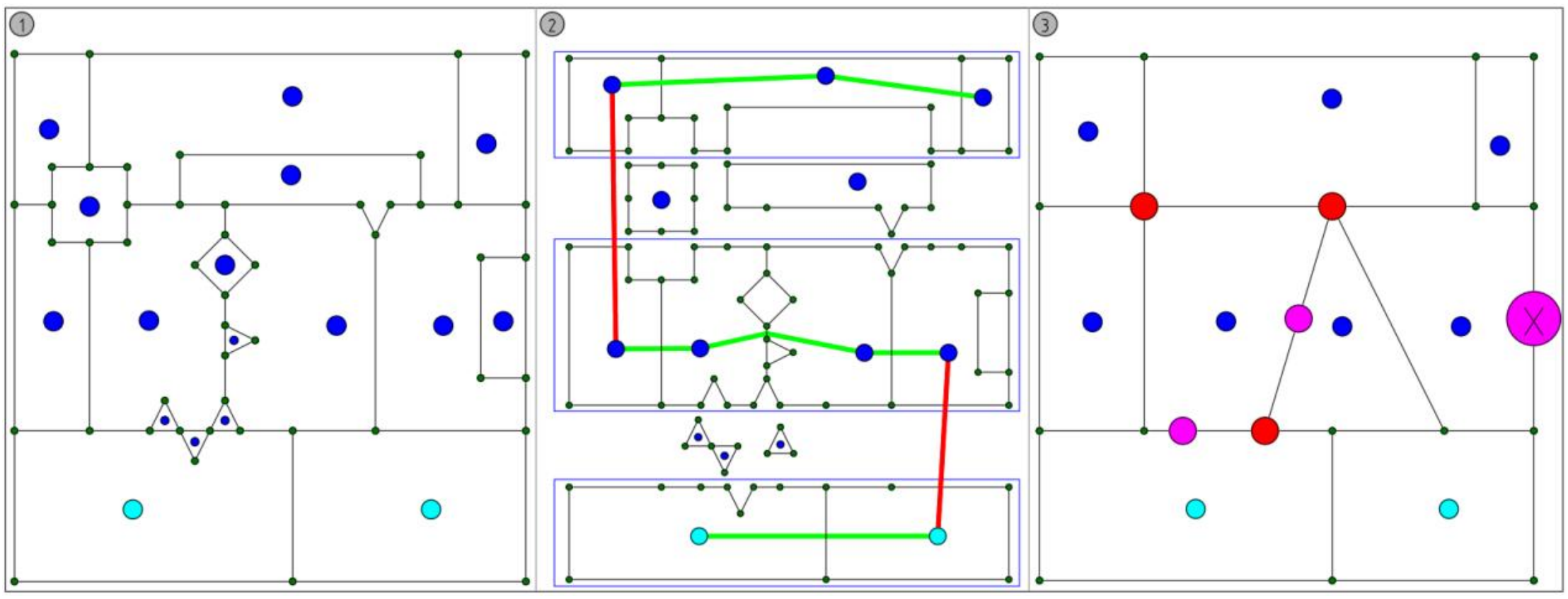

Fig.25 Strip division and simplified expression of an E0-graph

### 2 Main-route, C4-vertex and C4-chain

Main-route: As shown in Figure 25.2, on the E0-graph, the SPF calculation is supplemented after connection calculation to simplify the new strip. And the remaining room chain is the main-route for the new strip. The rooms on the chain are "main rooms", and adjacent "main rooms" have common edges. For the simplified rooms, according to the connection rule, they are automatically combined into C4s (inward-subgraphs).

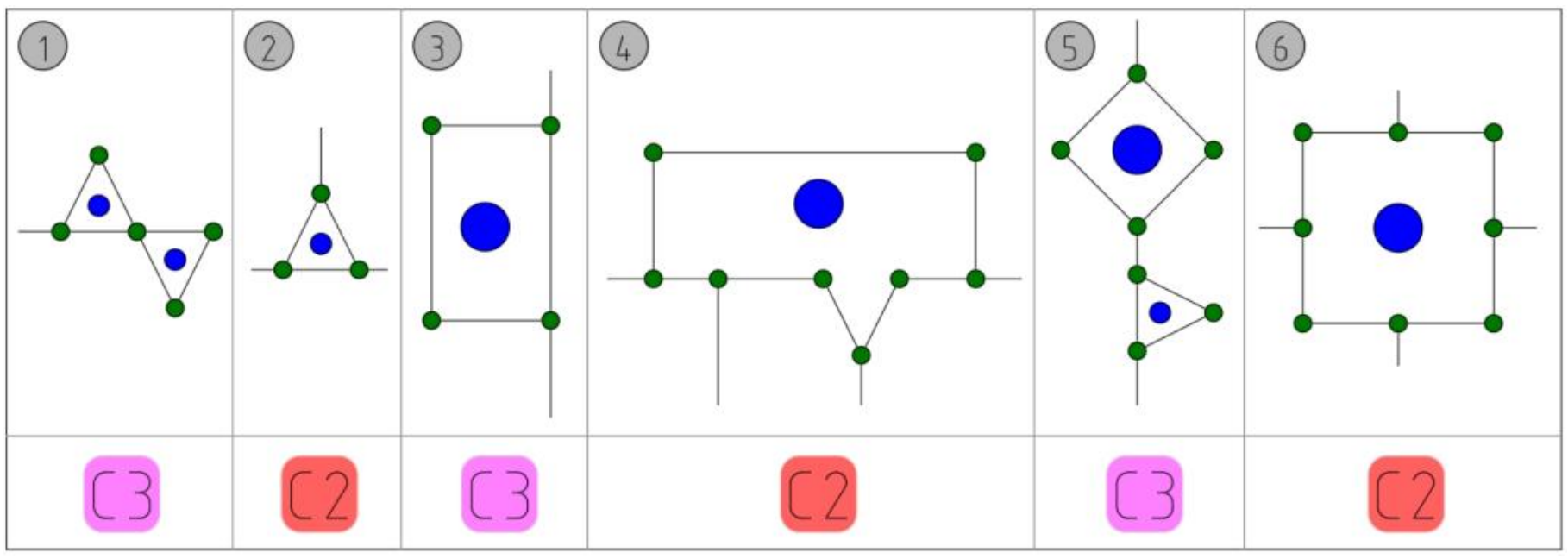

Fig.26 C4-vertices and C4-chains of an E0-graph

C4-vertex (C2) and C4-chain (C3): After the "main route" calculation, a C4-vertex refers to a C4 object serves as a common node of adjacent strips in the combined-grid; while a C4-chain refers to a group of C4 between a pair of adjacent main rooms, and they can be linked into a chain by the public edge of the two main rooms. The C3s and C4s in the E0-graph case are shown in Figure 26.

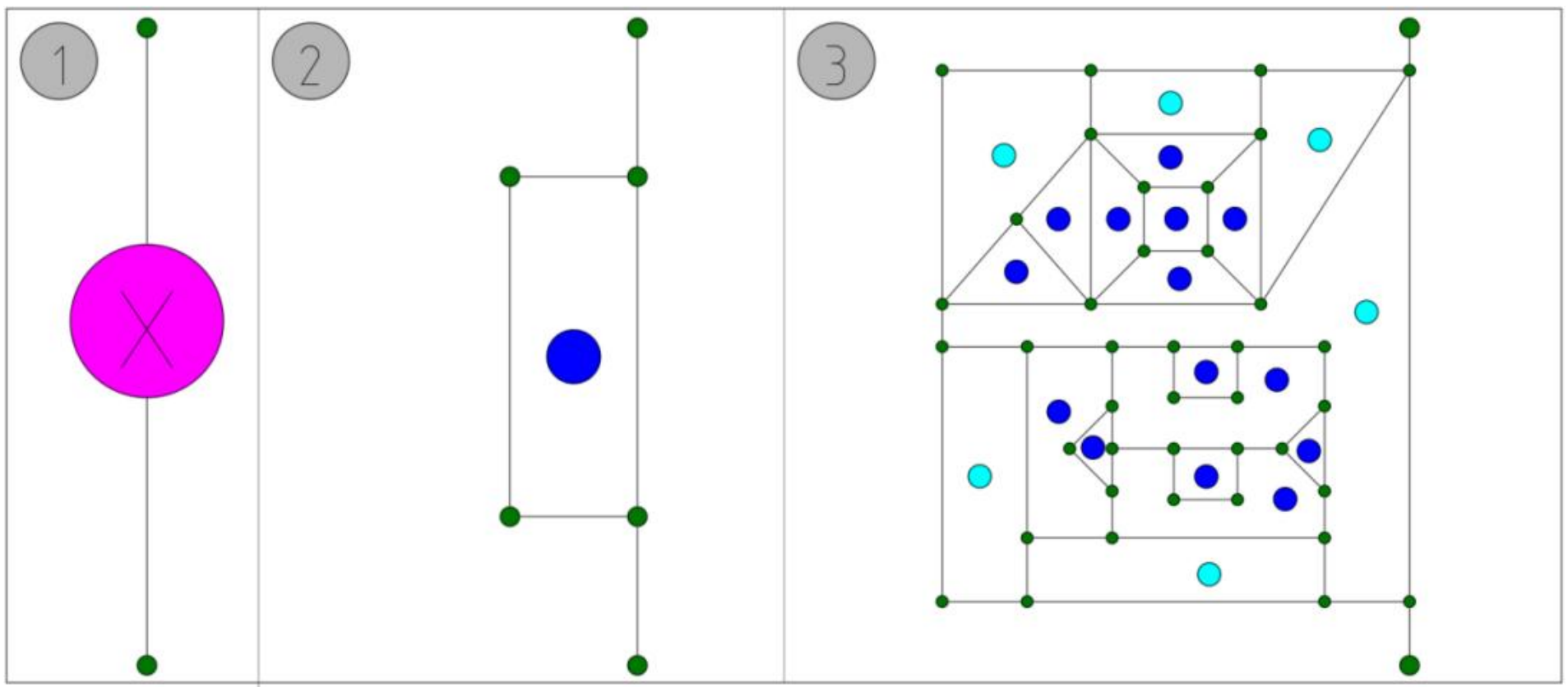

Fig.27 Refinement of a C4-chain

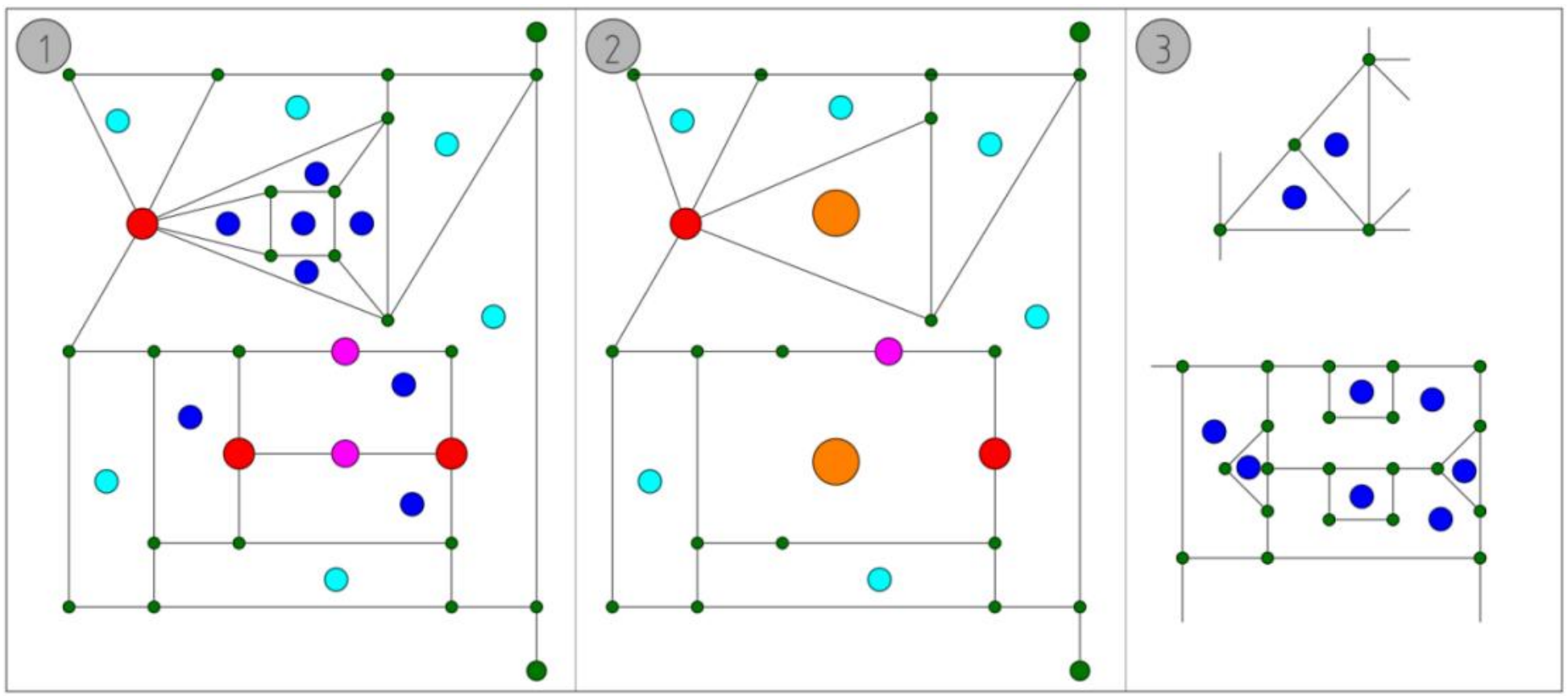

Fig.28 Composition of a refined C4-chain

Inward-main-route: The expansion and traversal of C3/C4 also need to perform the "main route" calculation, and its outer boundary is set "non-degenerate", that is, the vertex simplification of the outer boundary is not performed during the first-strip calculation. This setting can limit the complexity of the result’s structure when calculating the boundaries of relevant space objects. The refinement of the C4- chain and the calculation of its inward-main-route are shown in Figure 28.

Obviously, based on the main-route calculation and the concepts of C4-vertex and C4-chain, a E0-graph can be transformed into an expandable combined-grid.

**3 Syntax extensions and case expressions for 1-d graph**

For E0-graph, the 1-d graph system adds the following 9 symbols:

| S/N | symbol | Explanation |
|---|---|---|
| 1 | C2 | C4 as a vertex of the main-route |
| 2 | C3 | C4 chain as an edge of the main-route |
| 3 | A1 | public position-vertex of a C2 and a main-room in the same strip |
| 4 | A2 | public position-vertex of the adjacent main-room pair from a C3 |
| 5 | A3 | public position-vertex of a C2 and a secondary main-room |
| 6 | F5 | extensible F0 |
| 7 | F6 | extensible F1 |
| 8 | F7 | Identical edge boundary matching |
| 9 | F8 | extensible F2 |

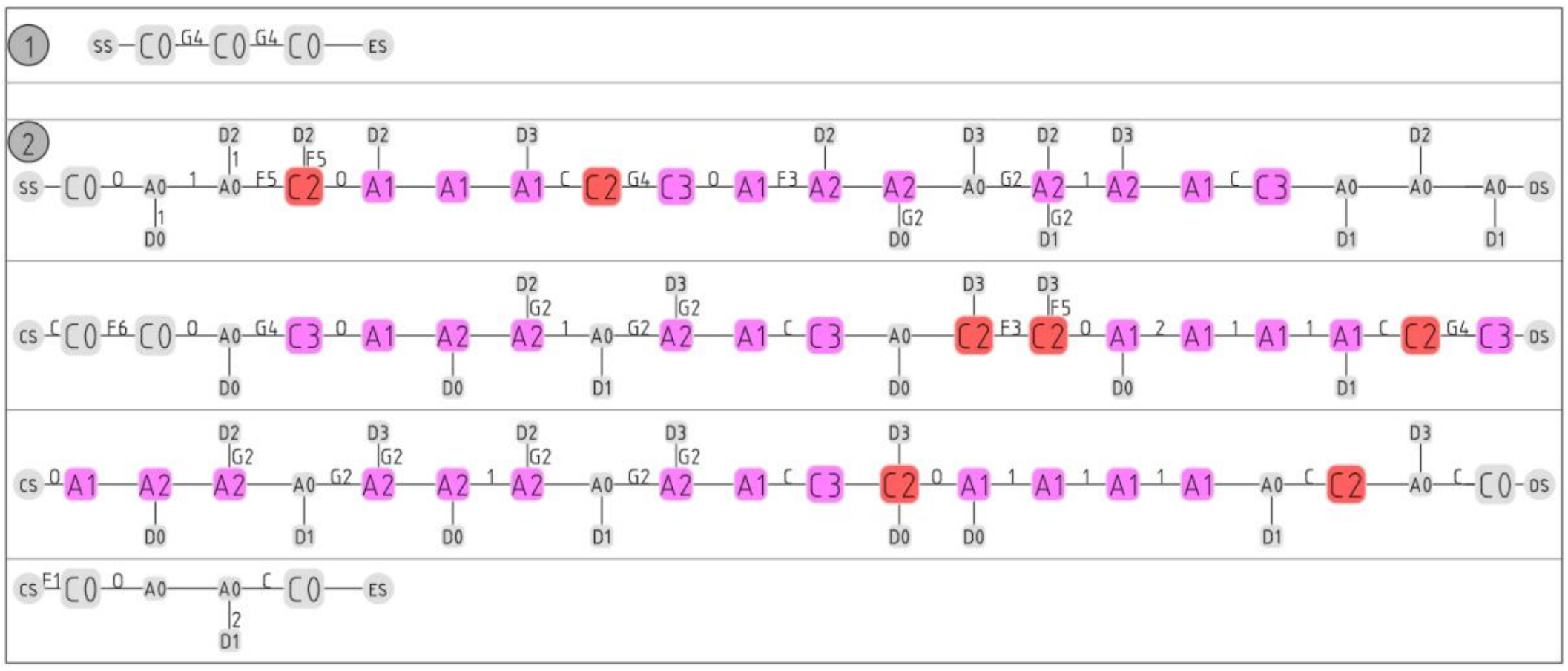

Fig.29 1-d representation of E0-graph

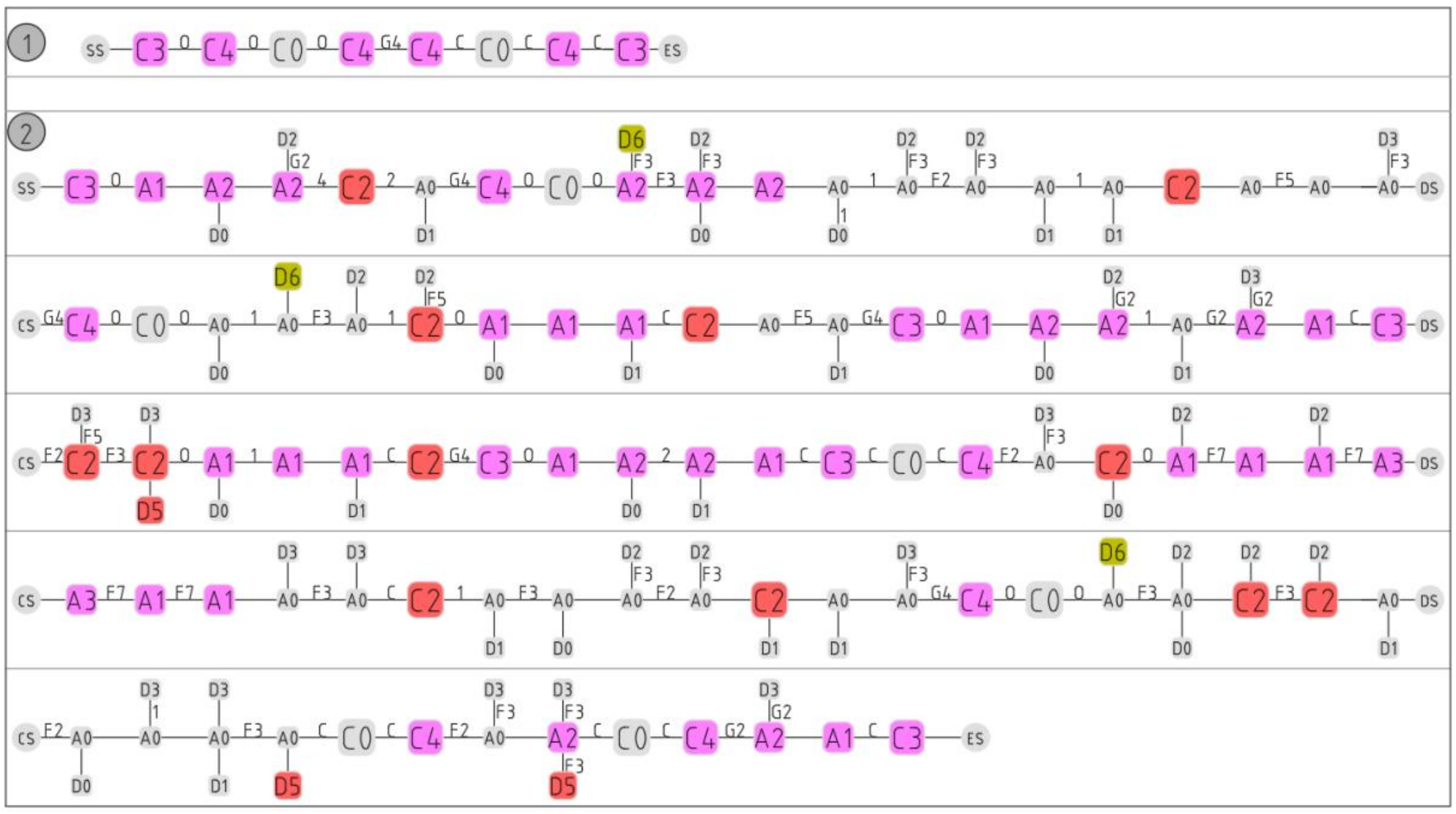

Fig.30 1-d representation of refined C4-chain

As shown in Figure 29.1 and Figure 30.1, the overall space composition is very simple, while the unfolded 1-d representation shown in Figure 29.2 and Figure 30.2 is relatively complex.

# 3 En-graph

**Chapter introduction**

There are 2 sections in this chapter. Section 1 completes the description of a En-graph’s composition by introducing the concepts and cases of embedding-level, combined-tree and tree-chained room. Section 2 demonstrates the 1-d representation for each concrete case.

## 3.1 Concept and Case

**1 "Simplify" Operation and En-graph**

The concept of "room with embedded structures" comes from a trivial fact that any planar graph can fit into a "large" room in a "composition" perspective.

For the embedded phenomenon, a "simplify" operation can be defined: for a given planar graph, delete all embedded vertices and edges until any vertex is not inside any room.

Furthermore, based on "Simplify" operation, the "En-graph" can be defined: a planar graph that can be determined as a E0-graph after the "simplify" operation.

**2 Embedding-level**

Embedding-level of room and En-graph: The "embedding-level" concept is similar to the "depth" concept of a tree. Referring to the "simplify" operation, for a given En-graph, the rooms in the simplified E0-graph are called "E0-rooms" of the original En-graph. Furthermore, for each subgraph in an embedded room, its remaining rooms after "simplify" calculation are called "E1-rooms". The calculation for all rooms’ E-value can be done by repeating the "simplify" operation. And the embedding value of a En-graph is set to be the maximum E-value of the room. (If the highest is 5, which means that there are at most 5 levels of embedded structures, the En-graph can be called a E5-graph.)

Based on the concept of "embedding-level", a parameter pair (EA and EB) can be defined: For any room in a given En-graph, its EA-value is the room’s E-value relative to the En-graph, and its EB-value is the E-value of the room itself as an subgraph. Furthermore, for any subgraph of the En-graph, its EA-value is the minimum EA-value of all its rooms, and its EB-value is the maximum EB-value of all its rooms.

Obviously, for any En-graph, if 2 rooms have common edges, their EA-values are equal. Based on this trivial observation, a subset of En-graph ("Set of En-C4") can be defined: For a certain room in a given En-graph, all other rooms in the graph with the same EA-value as this room can be extracted, and then the rooms would combine and divide themselves based on the requirement of "room continuity". Finally, a C4 set with the same EA-value can be obtained. (If the EA-value is 2, the set is called "E2-C4 set" of the En-graph)

**3 combined-tree and tree-chain room**

Combined-tree: In the S-rule system, the core role of "combined-tree" is as a logical tool for standardized expression between the rooms and its embedded C4s. Inside the room, a combined-tree can basically correspond to a continuous subgraph, but needs to be reconstructed into a tree structure. Taking the room as the En-graph, first extract its E1-C4 set, and then divide the set based on connection rule. Each C4 or discrete tree is treated as a vertex, and the edge between tree-vertex and C4-vertex is derived from their public position-vertex. For a public position-vertex of adjacent C4s, it will also derive a new tree-vertex if its degree >2 and and it doesn’t belong to any other tree-vertex. Finally, the entire embedded structures of a room are organized into multiple discrete combined-trees, some of which may have public positioning vertices with the outer boundary of the room.

EA-value and EB-value of a combined-tree: Suppose the EA-value of its belonging room is X, the combined-tree’s EA-value is X + 1, and its EB-value is the maximum EB-value of all the C4s it contains.

Tree-chained room: Emphasizes a simplified representation of embedded room, which means that a room expands into a combined-tree sequence. Furthermore, the "boundary" of a tree-chained room can also be correspondingly expressed as a combination of its outer boundary and the boundaries of its combined-trees. The boundary of a combined-tree is equivalent to the sequence obtained by the initial traversal of DFS order (meaning that the internal traversal of C4-vertices are not performed).

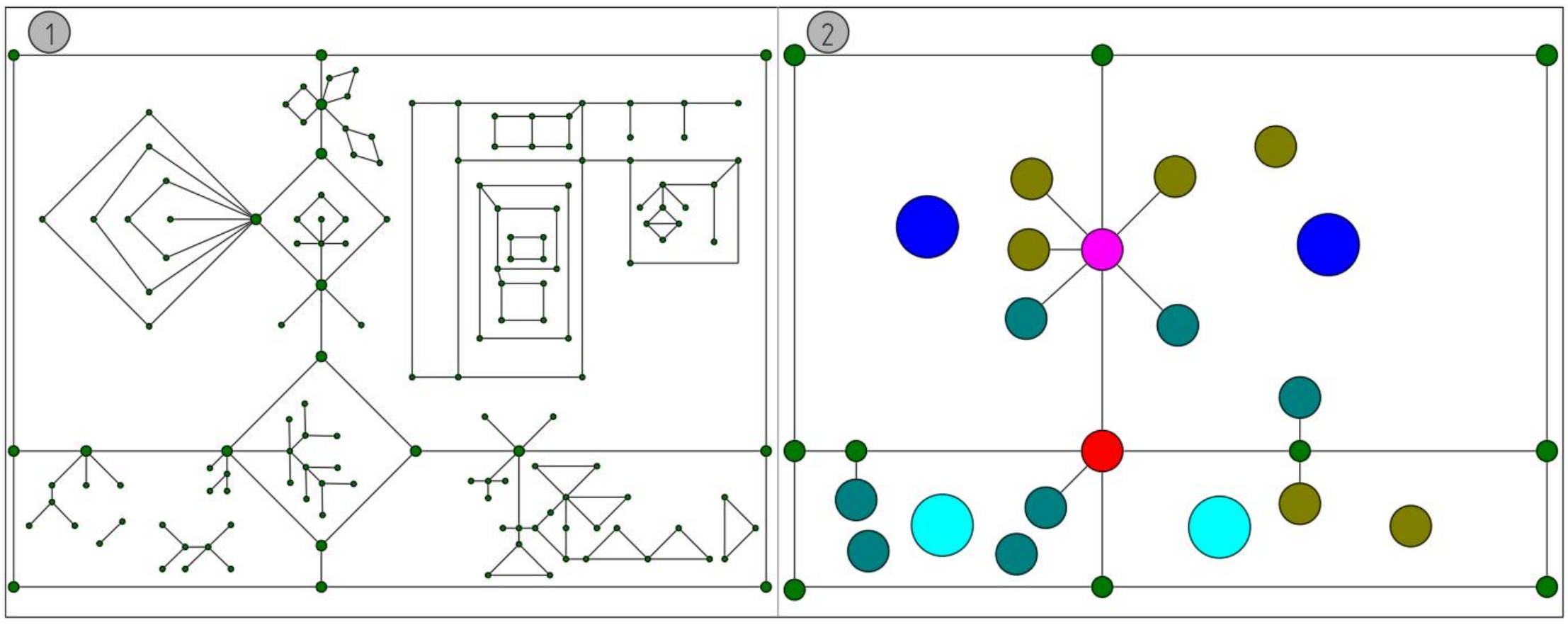

Fig.31 En-graph case and its simplified expression

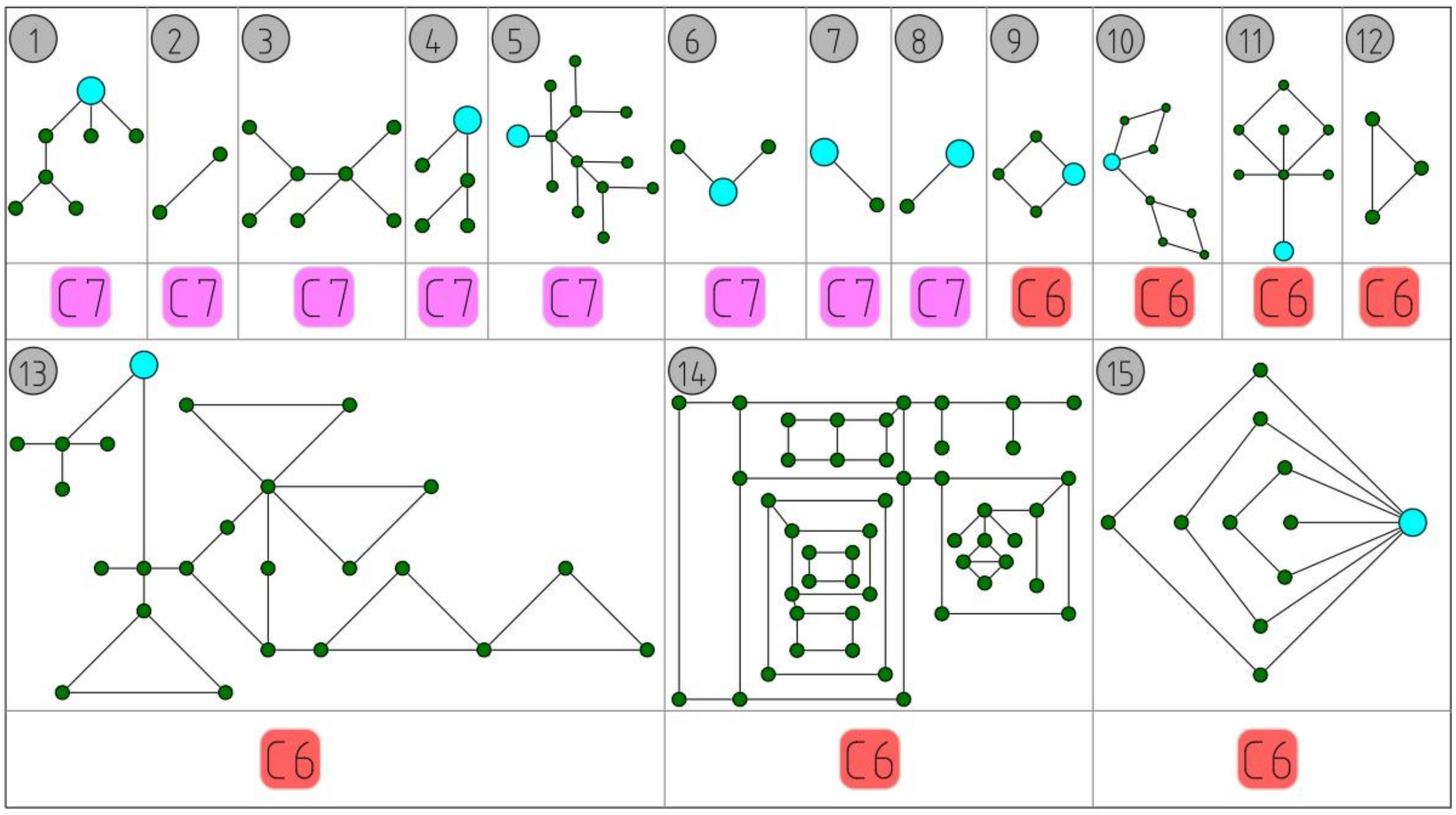

Fig.32 Combined-trees with EA value = 1 in the En-graph

The case of En-graph is shown in Figure 31.1, the corresponding E0-graph and combined-trees labels are shown in Figure 31.2. Taking the two rooms below as the initial data, the E0-C4 has 2 strips, including 4 main-rooms, 1 C2 and 1 C3. As shown in Figure 32, there are 15 combined-trees with EA value = 1, of which there are 13 in the main-rooms (7 of which contain only 1 tree-vertex). And there is 1 combined-tree in each of the C2 and C3.

**4 Composition and Expansion of En-graph**

There are 3 trivial observations about the composition of En-graph:

1 An En-graph can be regarded as an E0-graph in which the rooms may contain embedded subgraphs;

2 Each "embedded room" can be regarded as a combination of a standard room and a limited number of combined-trees with the same EA-value;

3 Each combined-tree can be regarded as a tree structure composed of tree-vertices and secondary En-graphs (C4-vertices).

On the basis of the above cognition, it can be found that when traversing an En-graph, its space-expansion will repeats in a fixed logical order: "tree-chained room" "combined-tree" "C4-vertex" "secondary tree-chained room"..., until all the embedded subgraps have been completely traversed.

## 3.2 1-d Expression of En-graph and Related Concepts

**1 Added symbols**

For En-graph, the 1-d graph system adds the following 4 symbols:

| S/N | Symbol | Explanation |
|---|---|---|
| 1 | C5 | "tree-chained room", room with embedded C6s |
| 2 | C6 | "combined tree", subgraph composed of trees (C7-vertices) and secondary En-graphs (C4-vertices) in a tree structure |
| 3 | G7 | Edge for C6, originates from the common positioning vertices of C4 and C7, the essence is an identical edge |
| 4 | G9 | The derived edge between a C5 and its C6s |

**2 combined-tree**

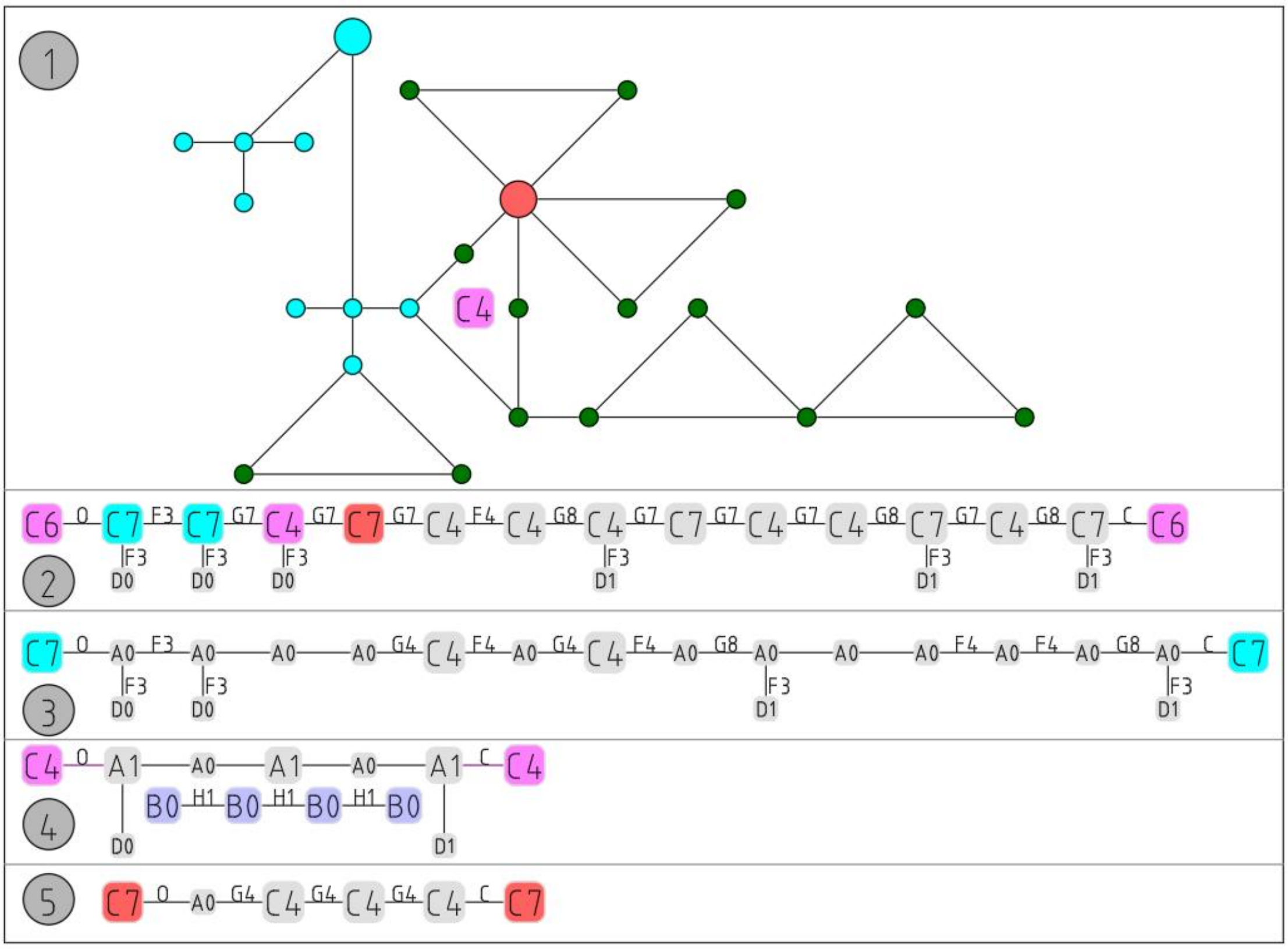

Fig.33 Traversal and 1-d representation of No.13 combined-tree

The traversal and 1-d expression of No. 13 combined-tree is shown in Figure 33:

Figure 33.2: Obviously, since the upper position-vertex is a common vertex with the room's outer boundary, the tree (C7-vertex) to which the common vertex belongs is the start of traversal.

Figure 33.3: In the expansion of the starting tree, the external related C4-vertices are echoed by interpolating the pairs of G4 and C4.

Figure 33.4: There are 3 A1-vertices on the boundary of the purple C4-vertex.

Figure 33.5: The red C7-vertex is derived from the public position-vertex of 3 C4-vertices.

## 3 Typical "embedding"

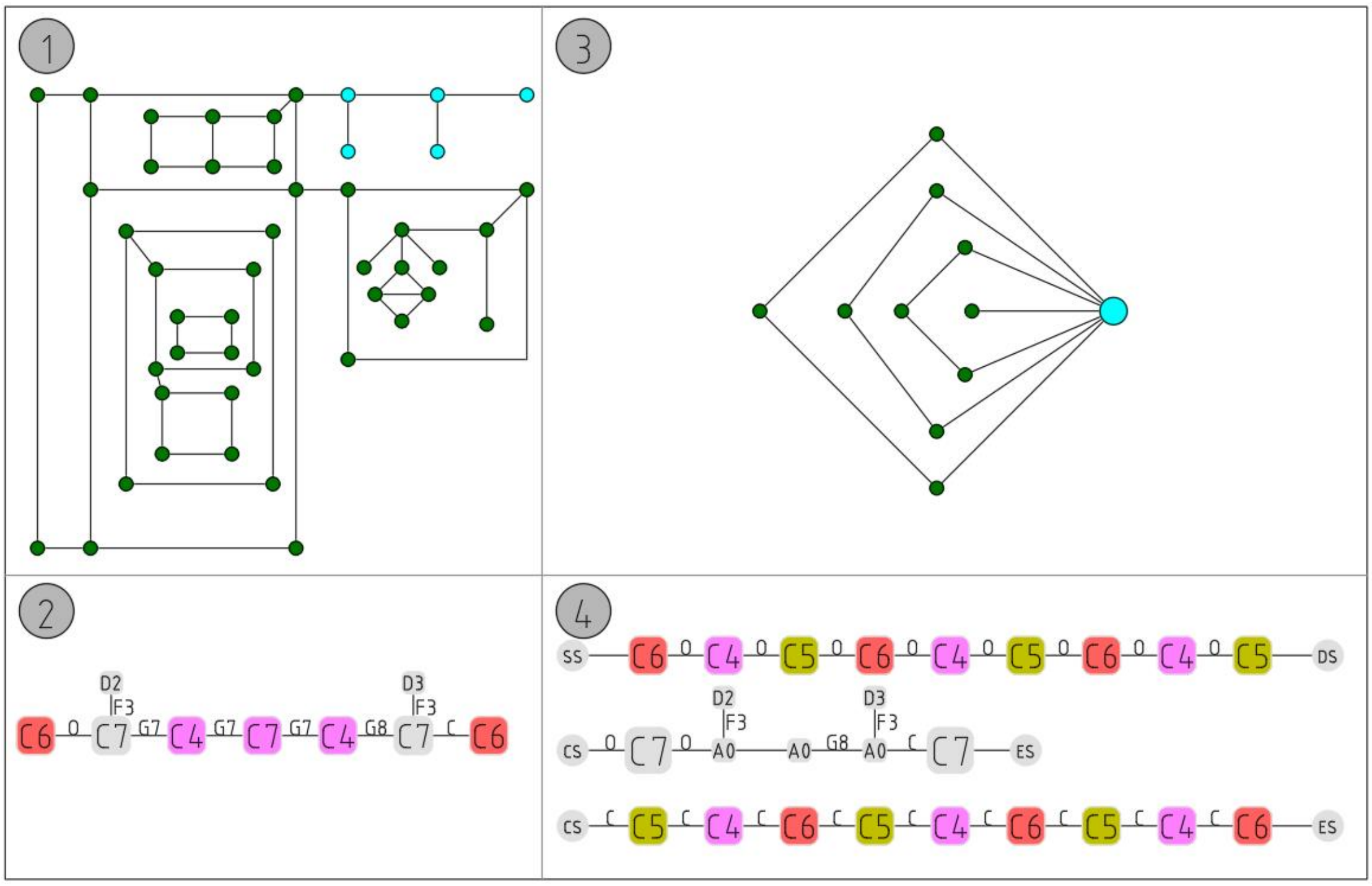

Fig.34 Combined-trees with EB value > 0 in the En-graph

The preliminary traversal and 1-d expression of No. 14-15 combined-trees are shown in Figure 34:

Figure 34.2: Corresponding to Figure 34.1, the No. 14 combined-tree is a chain composed of 2 C7-vertices and 2 C4-vertices.

Figure 33.4: Corresponding to Figure 34.3, No. 15 combined-tree has 3 embedding-levels, and the fixed logical sequence repeated 3 times can be seen in the preliminary 1-d representation, which means the "C6-C4-C5" sequence in the opening process and the "C5-C4-C6" sequence in the closing process.

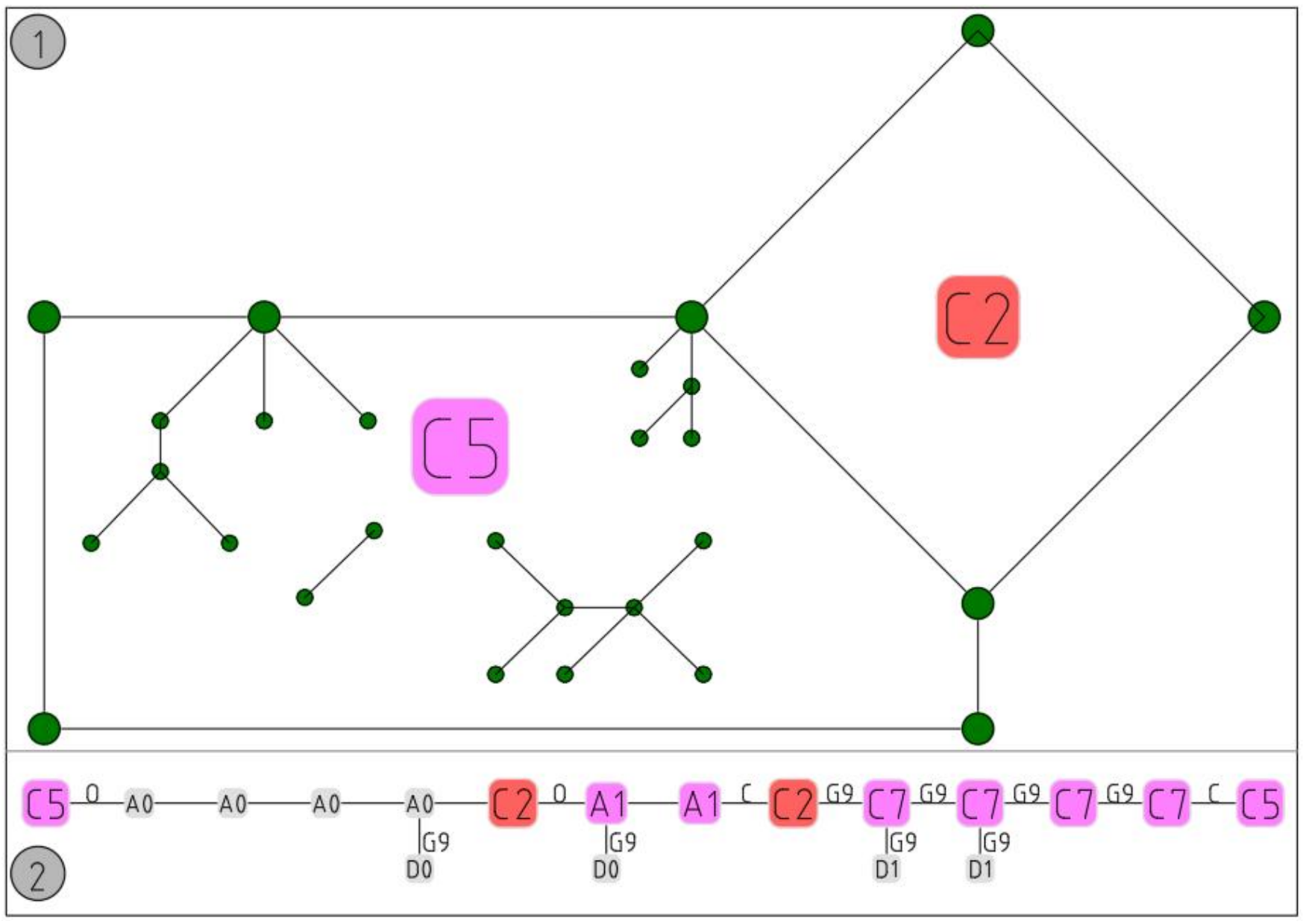

Fig.35 Main tree-chained room in the En-graph

The preliminary traversal and 1-d expression of the first main-room is shown in Figure 35. Note two syntax details:

1 Considering the matching rule of call-back pair, when arranging the C6-vertex chain, the start should be the end-vertex of the outer boundary, and C6-vertices with common position-vertices are preferentially arranged in the reversed direction.

2 When a common position-vertex is located in a C2 or C3, the C2 or C3 needs to be partially expanded to show it.

**4 Complete 1-d representation of En-graph**

Figures 36 and 37 show the 1-d representation of the En-graph case:

Figure 36.1 shows the overall structure, which is divided into 2 strips.

Figure 36.2 shows the preliminary representation in the form of an E0-graph.

Figure 37 shows the complete 1-d representation.

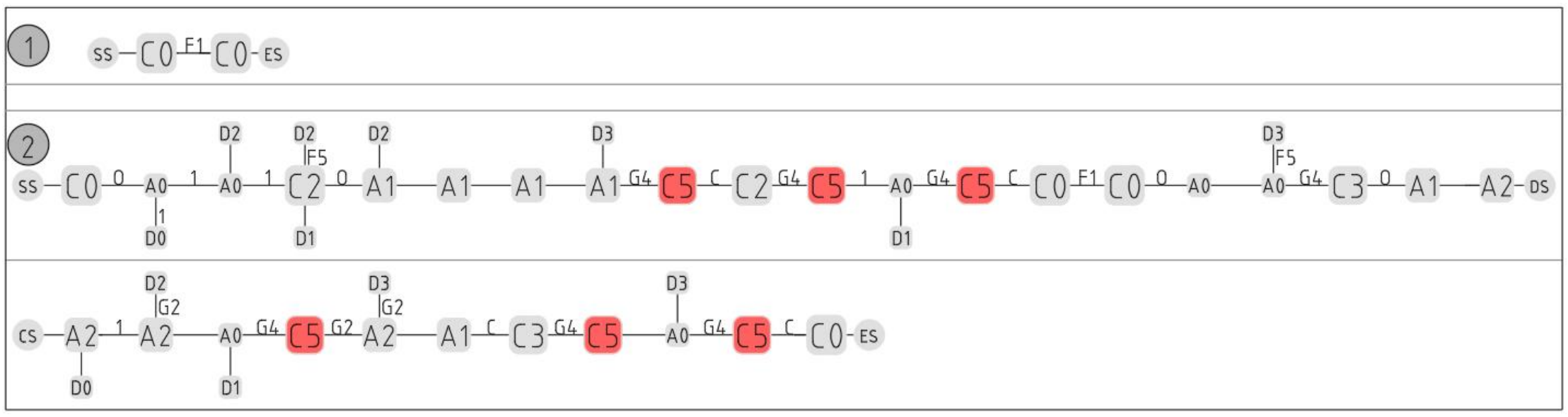

Fig.36 Preliminary 1-d representation of the En-graph

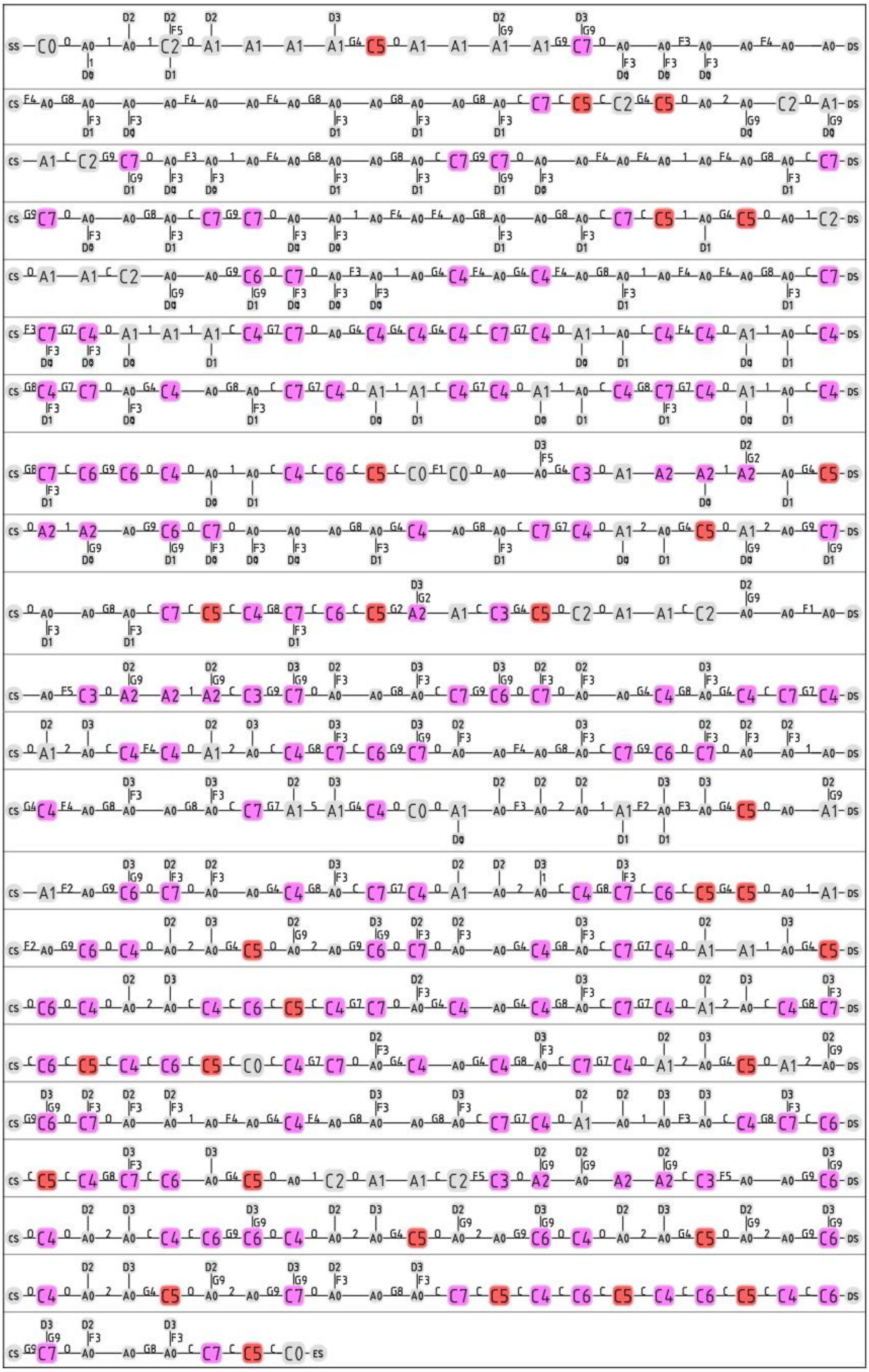

Fig.37 Complete 1-d representation of the En-graph

# 4 Summary

**Chapter introduction**

This chapter consists of three sections, aiming at summarizing the S-rule, main concepts and 1-d graph grammar established from grid to En-graph:

From the perspective of algorithms, the S-rule is an enhanced DFS, and the 1-d graph is a standard format for data recording in the process of "reading" a planar graph with S-rule.

After satisfying the requirement of "index-free", the 1-d graph gradually conforms to the typical cognition of "language" by AEC personnel and algorithms: based on the simple combination of finite types of symbols, a standardized sequence is given in a unified direction. In turn, it may become the basis for efficient interactive manipulation of planar graph data.

## 4.1 Summary of S-Rule

**1 The dual nature of S-rule**

As the basis of traversal and cognition, a simple and "objective" division of the planar graph is required first. The primary nature of "S-rule" is a rule that "objectively" divides the planar graph based on initial data, that is, the grid is used as the ideal format, and the "strip" concept (C0) is established based on its rows/columns. And the core assumption is that each planar graph with "room continuity" can be expressed as a continuous strip sequence.

The secondary nature of "S-rule" is the traversal order for this division, and the core principle of establishing the "traversal order" is to follow the AEC logic and the original intention of the "language". That is, the direction of the continuous room sequence under the "strip sequence" framework is used as the "guidance" for the direction of the derived vertex sequence.

**2 The "objectivity" of S-Rule**

From the perspective of AEC, the row division of a grid is a division rule that derives the "average" strip sequence based on the "average" initial data. And for more general planar graph and initial data, "averageness" translates to "objectivity".

That is, for all kinds of secondary space objects derived during the traversal process, the S-rule uniformly follows the "inward" principle, so as to avoid presetting specific secondary initial data for them. Therefore, the excessive influence of the initial data on the final division result is limited.

## 4.2 The 2 core concept pairs serving the S-rule

Except for the strip concept, other space concepts (Cx) are derived from the non-standard nature of the planar graph relative to the grid, and their purpose is still to serve the strip sequence. For E0-graph and En-graph, each has a core concept pair.

**1 "inward-subgraph" and "tree" serving the E0-graph**

Based on the space concept system, the S-rule can transform a E0-graph into a specific type of combined grid, in which vertices, edges and rooms can be expanded.

Inward-subgraph is the public definition serving C4-vertex, C4-chain, and general inward combined-grid, and it is also a unified option that can serve the "objectivity" requirement of S-rules.

Tree with DFS are the underlying logic for combining the main-route and the concept system of inward-subgraph: For a general E0-graph (which can be regarded as an C4 as a whole), when it is transformed into an expandable combined-grid after the main-route calculation, the final result is usually a tree structure consisting of strips (C0) and branches (C1) or secondary inward-subgraphs (C4), and the expansion of these objects would still follow the identical rule.

**2 "combined-tree" and "tree-chained room" serving the En-graph**

Obviously, regardless of the embedded structure, an E0-graph can be regarded as a C4-vertex of a combined-tree inside a certain room, and can also be regarded as an object in a set of En-C4.

After the preliminary traversal of the C4 is completed, S-rule use two combined methods to complete the division and traversal of the embedded rooms and the upper embedded-level subgraphs. For an embedded rooms, first complete the traversal from the room to its combined-trees with "downward" rule, and then complete the traversal of each C4 with "inward" rule. For a traversed combined-tree, first complete the traversal of the upper embedded-level room with "upward" rule, and then complete the traversal of the upper-level subgraph with "outward" rule. The two processes will be repeated alternately until the full graph traversal is completed.

The core function of the "concept pair" is to ensure that the 1-d expression maintains a fixed symbolic logic order when it crosses the embedded boundary multiple times in the "downward" or "upward" mode.

## 4.3 Core rule and 4 basic grammars of 1-d graph

"Connect-forward" as the core rule: There is a connection relationship or a simplified connection relationship that has been traversed between adjacent position vertices or space vertices in the 1-d graph. In each strip, the edge direction keeps forward according to the direction of the rooms, that is, it always points to the next untraversed room. On this basis, the 1-d graph sets 4 basic grammars:

1 call-back: The outer boundary of each room, i.e. its terminating edge on the 1-d graph, is marked by a call-back pair. Furthermore, when traversing between space objects, in order to ensure the continuity of traversal, the backtracking matching of the common vertex is also ensured by a call-back pair.

2 Derive-traversal: each derived space vertex is expanded and internally traversed at the first traversal position.

3 open-close: The internal traversal of each space vertex is marked by a pair of open and close edges.

4 Embed-arrangement: When traversing an embedded room "inwardly", the 1-d graph first arranges its combined-trees as a chain on the right; When traversing the parent embedded room of a combined-tree "outwardly", the 1-d graph first derives its outer boundary on the left, and then arranges the sibling combined-trees as a chain on the right.